\documentclass[10pt,a4paper]{article}

\usepackage{makeidx}
\usepackage{latexsym}
\usepackage{amsfonts}
\usepackage{amsmath}
\usepackage{amsthm}
\usepackage{amssymb}
\usepackage{amscd}
\usepackage[dvips]{graphicx}

\theoremstyle{definition}
\newtheorem{theorem}{Theorem}[section]
\newtheorem{prop}[theorem]{Proposition}
\newtheorem{lemma}[theorem]{Lemma}
\newtheorem{corollary}[theorem]{Corollary}

\newtheorem{remark}[theorem]{Remark}

\newenvironment{demo}[1]{%
  \trivlist
  \item[\hskip\labelsep
        {\bf #1.}]
}{%
\hfill\qedsymbol
  \endtrivlist
}

\newenvironment{roster}[1]{%
\begin{list}{}{%
\setlength{\topsep}{0pt}
\setlength{\itemsep}{0pt}
\setlength{\parsep}{0pt}
\settowidth{\labelwidth}{#1}
\setlength{\leftmargin}{\labelwidth}
\addtolength{\leftmargin}{\parindent}
}}{\end{list}}
\newcommand\Int{\mathbb{Z}}
\newcommand\Nat{\mathbb{N}}

\renewcommand\tilde{\widetilde}
\newcommand\Pf{\operatorname{Pf}}

\newcommand\Sym{{\mathcal{S}}}

\newcommand\sgn{\operatorname{sgn}}
\newcommand{\Par}{{\mathcal{P}}}
\newcommand{\Q}{{\mathcal{Q}}}
\newcommand{\R}{{\mathcal{R}}}
\newcommand\rdots{\mathinner{\mkern1mu\raise0pt\vbox{\kern7pt\hbox{.}}
     \mkern2mu\raise4pt\hbox{.}\mkern2mu\raise8pt\hbox{.}\mkern1mu}}
\newcommand\LR{\operatorname{LR}}     
\newcommand\trans{{}^t\!}             

\newcommand{\vectx}{\boldsymbol{x}}
\newcommand{\vecty}{\boldsymbol{y}}
\newcommand{\vectz}{\boldsymbol{z}}
\newcommand{\vectw}{\boldsymbol{w}}
\newcommand{\vecta}{\boldsymbol{a}}
\newcommand{\vectb}{\boldsymbol{b}}
\newcommand{\vectc}{\boldsymbol{c}}
\newcommand{\vectd}{\boldsymbol{d}}
\newcommand{\vectxi}{\boldsymbol{\xi}}
\newcommand{\vecteta}{\boldsymbol{\eta}}
\newcommand{\vectzeta}{\boldsymbol{\zeta}}
\newcommand{\vectomega}{\boldsymbol{\omega}}
\newcommand{\vectalpha}{\boldsymbol{\alpha}}
\newcommand{\vectbeta}{\boldsymbol{\beta}}
\newcommand{\vectgamma}{\boldsymbol{\gamma}}
\newcommand{\vectdelta}{\boldsymbol{\delta}}
\newcommand{\vectone}{\boldsymbol{1}}

\newcommand{\vectcp}{\boldsymbol{c'}}

\def\U#1#2#3#4#5#6{%
U^{{#1},{#2}} \left( \hspace{-\arraycolsep}
 \begin{array}{c|c}
  {#3} & {#5} \\ {#4} & {#6}
 \end{array}
 \hspace{-\arraycolsep} \right)}

\numberwithin{equation}{section}
\title{
Generalizations of Cauchy's Determinant and Schur's Pfaffian
}

\author{
Masao Ishikawa\thanks{
Faculty of Education,
Tottori University,
 e-mail:\tt{ishikawa@fed.tottori-u.ac.jp}}
\and
Soichi Okada\thanks{Graduate School of Mathematics, Nagoya University,
 e-mail:\tt{okada@math.nagoya-u.ac.jp}}
\and
Hiroyuki Tagawa\thanks{Faculty of Education, Wakayama University,
 e-mail:\tt{tagawa@math.edu.wakayama-u.ac.jp}} 
\and
Jiang Zeng\thanks{Institut Girard Desargues, Universit\'e Claude Bernard Lyon-I,
 e-mail:\tt{zeng@igd.univ-lyon1.fr}} 
}

\date{\empty}

\begin{document}

\maketitle

\begin{abstract}
We present several identities of Cauchy-type determinants
 and Schur-type Pfaffians involving generalized Vandermonde
 determinants, which generalize Cauchy's determinant
 $\det \left( 1/(x_i+y_j) \right)$ and
 Schur's Pfaffian $\Pf \left( (x_j - x_i)/(x_j + x_i) \right)$.
Some special cases of these identities are given by
 S.~Okada and T.~Sundquist.
As an application, we give a relation for the Littlewood--Richardson
 coefficients involving a rectangular partition.
\end{abstract}

\bigbreak

\noindent{\bf Keywords} : Pfaffian,  Cauchy's Determinant,  Pl\"ucker relations,  Schur functions.


\section{Introduction}


Computations of  determinants and Pfaffians are of great
importance not only in many branches of mathematics but also in
physics. Some people need relations among minors or subpfaffians
of a general
 matrix, others have to evaluate special determinants or Pfaffians.
In enumerative combinatorics and representation theory, a central
role is played by Cauchy's determinant identity~\cite{C}
\begin{equation}
\det \left(
 \frac{ 1 }{ x_i + y_j }
\right)_{1 \le i, j \le n}
 =
\frac{ \prod_{1 \le i < j \le n} (x_j - x_i) (y_j - y_i) }
     { \prod_{i,j=1}^n (x_i + y_j) },
\label{cauchy}
\end{equation}
and Schur's Pfaffian identity \cite{S}
\begin{equation}
\Pf \left( \frac{x_j - x_i}{x_j + x_i} \right)_{1 \le i, j \le 2n}
 =
\prod_{1 \le i < j \le 2n} \frac{x_j - x_i}{x_j+x_i}.
\label{schur}
\end{equation}
The reader is referred to \cite{I,IW2,LLT,O1,O2,Ste,Su} for some
recent  variations and generalizations with their applications of
\eqref{schur} and \eqref{cauchy}. Besides, Krattenthaler~\cite{Kr}
has given a comprehensive survey of determinant evaluations.


In the same vein, we shall give several identities of Cauchy-type
determinants and
 Schur-type Pfaffians whose entries involve two kinds of  generalized Vandermonde determinants.
Let $\vectx = (x_1, \cdots, x_n)$ and $\vecta = (a_1, \cdots, a_n)$
 be two vectors of variables of length $n$. Let $p$ and $q$ be
two nonnegative integers  such that $p+q = n$.
 Denote by  $V^{p,q}(\vectx ;
 \vecta)$  the $n \times n$ matrix with $i$th row
$$
(1, x_i, \cdots, x_i^{p-1}, a_i, a_i x_i, \cdots, a_i x_i^{q-1}),
$$
and $W^n(\vectx ; \vecta)$ the $n \times n$ matrix with $i$th row
$$
(1+a_i x_i^{n-1}, x_i + a_i x_i^{n-2}, \cdots, x_i^{n-1} + a_i).
$$
For example, if $q=0$, then $V^{n,0}(\vectx ; \vecta)
 = \left( x_i^{j-1} \right)_{1 \le i, j \le n}$ is the usual
 Vandermonde matrix and $\det V^{n,0}(\vectx ; \vecta)
 = \prod_{1 \le i < j \le n} (x_j - x_i)$.
If $p=q=1$, then $\det V^{1,1}(\vectx; \vecta) = a_2 - a_1$, while
the matrices $V^{3,2}(\vectx;\vecta)$ and $W^5(\vectx;\vecta)$
 can be visualized as follows:
\begin{gather*}
V^{3,2}(\vectx ; \vecta)
 =
\begin{pmatrix}
 1 & x_{1} & x_{1}^{2} & a_{1} & a_{1}x_{1} \\
 1 & x_{2} & x_{2}^{2} & a_{2} & a_{2}x_{2} \\
 1 & x_{3} & x_{3}^{2} & a_{3} & a_{3}x_{2} \\
 1 & x_{4} & x_{4}^{2} & a_{4} & a_{4}x_{4} \\
 1 & x_{5} & x_{5}^{2} & a_{5} & a_{5}x_{5}
\end{pmatrix},
\\
W^5(\vectx ; \vecta)
 =
\begin{pmatrix}
1 + a_1 x_1^4 & x_1 + a_1 x_1^3 & x_1^2 + a_1 x_1^2 & x_1^3 + a_1 x_1
 & x_1^4 + a_1 \\
1 + a_2 x_2^4 & x_2 + a_2 x_2^3 & x_2^2 + a_2 x_2^2 & x_2^3 + a_2 x_2
 & x_2^4 + a_2 \\
1 + a_3 x_3^4 & x_3 + a_3 x_3^3 & x_3^2 + a_3 x_3^2 & x_3^3 + a_3 x_3
 & x_3^4 + a_3 \\
1 + a_4 x_4^4 & x_4 + a_4 x_4^3 & x_4^2 + a_4 x_4^2 & x_4^3 + a_4 x_4
 & x_4^4 + a_4 \\
1 + a_5 x_5^4 & x_5 + a_5 x_5^3 & x_5^2 + a_5 x_5^2 & x_5^3 + a_5 x_5
 & x_5^4 + a_5
\end{pmatrix}.
\end{gather*}


The main purpose of this paper is to prove the following identities
 for the determinants and Pfaffians whose entries involve these generalized Vandermonde
 determinants.

\begin{theorem} \label{thm:main}
\begin{roster}{(d)}
\item[(a)]
Let $n$ be a positive integer and let $p$ and $q$ be nonnegative integers.
For six vectors of variables
\begin{gather*}
\vectx = (x_1, \cdots, x_n),\
\vecty = (y_1, \cdots, y_n),\
\vecta = (a_1, \cdots, a_n),\
\vectb = (b_1, \cdots, b_n),
\\
\vectz = (z_1, \cdots, z_{p+q}),\
\vectc = (c_1, \cdots, c_{p+q}),
\end{gather*}
we have
\begin{multline}
\det \left(
 \frac{ \det V^{p+1,q+1}(x_i,y_j,\vectz ; a_i,b_j,\vectc) }
      {y_j - x_i}
 \right)_{1 \le i, j \le n}
\\
=
\frac{ (-1)^{n(n-1)/2} }{ \prod_{i,j=1}^n (y_j - x_i) }
 \det V^{p,q}(\vectz ; \vectc)^{n-1}
 \det V^{n+p,n+q}(\vectx,\vecty,\vectz ; \vecta,\vectb,\vectc).
\label{main1}
\end{multline}
\item[(b)]
Let $n$ be a positive integer and
 let $p$, $q$, $r$, $s$ be nonnegative integers.
For seven vectors of variables
\begin{gather*}
\vectx = (x_1, \cdots, x_{2n}),\
\vecta = (a_1, \cdots, a_{2n}),\
\vectb = (b_1, \cdots, b_{2n}),
\\
\vectz = (z_1, \cdots, z_{p+q}),\
\vectc = (c_1, \cdots, c_{p+q}),
\\
\vectw = (w_1, \cdots, w_{r+s}),\
\vectd = (d_1, \cdots, d_{r+s}),
\end{gather*}
we have
\begin{multline}
\Pf \left(
 \frac{ \det V^{p+1,q+1}(x_i,x_j,\vectz ; a_i,a_j,\vectc)
        \det V^{r+1,s+1}(x_i,x_j,\vectw ; b_i,b_j,\vectd) }
      { x_j - x_i }
\right)_{1 \le i, j \le 2n}
\\
=
\frac{1}{\prod_{1 \le i < j \le 2n}(x_j - x_i)}
 \det V^{p,q}(\vectz ; \vectc)^{n-1}
 \det V^{r,s}(\vectw ; \vectd)^{n-1}
\\
\times
 \det V^{n+p,n+q}(\vectx,\vectz ; \vecta,\vectc)
 \det V^{n+r,n+s}(\vectx,\vectw ; \vectb,\vectd).
\label{main2}
\end{multline}
\item[(c)]
Let $n$ be a positive integer and let $p$ be a nonnegative integer.
For six vectors of variables
\begin{gather*}
\vectx = (x_1, \cdots, x_n),\
\vecty = (y_1, \cdots, y_n),\
\vecta = (a_1, \cdots, a_n),\
\vectb = (b_1, \cdots, b_n),
\\
\vectz = (z_1, \cdots, z_p),\
\vectc = (c_1, \cdots, c_p),
\end{gather*}
we have
\begin{multline}
\det \left(
 \frac{ \det W^{p+2}(x_i,y_j,\vectz ; a_i,b_j,\vectc) }
      { (y_j - x_i)(1-x_i y_j) }
 \right)_{1 \le i, j \le n}
\\
=
\frac{ 1 }{ \prod_{i,j=1}^n (y_j - x_i)(1-x_iy_j) }
 \det W^p(\vectz ; \vectc)^{n-1}
 \det W^{2n+p}(\vectx,\vecty,\vectz ; \vecta,\vectb,\vectc).
\label{main3}
\end{multline}
\item[(d)]
Let $n$ be a positive integer and let $p$ and $q$ be nonnegative integers.
For seven vectors of variables
\begin{gather*}
\vectx = (x_{1}, \cdots, x_{2n}),\
\vecta = (a_{1}, \cdots, a_{2n}),\
\vectb = (b_{1}, \cdots, b_{2n}),
\\
\vectz = (z_1, \cdots, z_p),\
\vectc = (c_1, \cdots, c_p),
\\
\vectw = (w_1, \cdots, w_q),\
\vectd = (d_1, \cdots, d_q),
\end{gather*}
we have
\begin{multline}
\Pf \left(
 \frac{ \det W^{p+2}(x_i,x_j,\vectz ; a_i,a_j,\vectc)
        \det W^{q+2}(x_i,x_j,\vectw ; b_i,b_j,\vectd) }
      { (x_j - x_i)(1 - x_i x_j) }
\right)_{1 \le i, j \le 2n}
\\
=
\frac{1}{\prod_{1 \le i < j \le 2n}(x_j - x_i)(1 - x_i x_j) }
 \det W^p(\vectz ; \vectc)^{n-1}
 \det W^q(\vectw ; \vectd)^{n-1}
\\
\times
 \det W^{2n+p}(\vectx,\vectz ; \vecta,\vectc)
 \det W^{2n+q}(\vectx,\vectw ; \vectb,\vectd).
\label{main4}
\end{multline}
\end{roster}
\end{theorem}

These identities were conjectured by one of the authors \cite{O3}.
If we put $p=q=0$ in (\ref{main1}) or $p=q=r=s=0$ in (\ref{main2}),
 then the identities read
\begin{gather}
\det \left( \frac{b_j - a_i}{y_j - x_i} \right)_{1 \le i, j \le n}
 =
\frac{ (-1)^{n(n-1)/2} }
     { \prod_{i,j=1}^n (y_j - x_i) }
\det V^{n,n}(\vectx,\vecty;\vecta,\vectb),
\label{special1}
\\
\Pf \left( \frac{(a_j-a_i)(b_j-b_i)}{x_j-x_i} \right)_{1 \le i,j \le 2n}
 =
\frac{ 1 }
     { \prod_{1 \le i < j \le 2n} (x_j - x_i) }
\det V^{n,n}(\vectx;\vecta) \det V^{n,n}(\vectx;\vectb).
\label{special2}
\end{gather}
These particular cases, as well as the identities (\ref{main3}) with $p=0$
 and (\ref{main4}) with $p=q=0$,
 are first given by S.~Okada \cite[Theorems~4.2, 4.7, 4.3, 4.4]{O1}
 in his study of rectangular-shaped representations of classical groups.
Another special case of the identity (\ref{main3}) with $p=1$
 is given in \cite{O2} and applied to the enumeration of
 vertically and horizontally symmetric alternating sign matrices.
These special cases are the starting point of our study.

Under the specialization
$$
x_i \leftarrow x_i^2,
\quad
y_i \leftarrow y_i^2,
\quad
z_i \leftarrow z_i^2,
\quad
w_i \leftarrow w_i^2,
\quad
a_i \leftarrow x_i,
\quad
b_i \leftarrow y_i,
\quad
c_i \leftarrow z_i,
\quad
d_i \leftarrow w_i,
$$
one can deduce from (\ref{main1}) and (\ref{main2}) the following identities:
\begin{align}
&
\det \left(
 \frac{ s_{\delta(k)}(x_i,y_j,\vectz) }{x_i + y_j}
\right)_{1 \le i, j \le n}
 =
 \frac{ \prod_{1 \le i < j \le n} (x_j - x_i) (y_j - y_i) }
      { \prod_{i,j=1}^n (x_i + y_j) }
 s_{\delta(k)}(\vectz)^{n-1} s_{\delta(k)}(\vectx,\vecty,\vectz),
\label{cauchy1}
\\
&
\Pf \left(
 \frac{ x_j - x_i }{ x_j + x_i }
 s_{\delta(k)}(x_i,x_j \vectz) s_{\delta(l)}(x_i,x_j,\vectw)
\right)_{1 \le i, j \le 2n}
\notag
\\
&\quad=
 \prod_{1 \le i < j \le 2n} \frac{ x_j - x_i }{ x_j + x_i }
 s_{\delta(k)}(\vectz)^{n-1} s_{\delta(l)}(\vectw)^{n-1}
 s_{\delta(k)}(\vectx,\vectz) s_{\delta(l)}(\vectx,\vectw),
\label{schur1}
\end{align}
where $s_\lambda$ denotes the Schur function corresponding to
 a partition $\lambda$ and $\delta(k) = (k,k-1,\cdots,1)$ denotes
 the staircase partition.
If we take $k=0$ in (\ref{cauchy1}) and $k=l=0$ in (\ref{schur1}),
 we obtain Cauchy's determinant identity (\ref{cauchy}) and
 Schur's Pfaffian identity (\ref{schur}).
Another special case of (\ref{cauchy1}) with $k=1$ is
 the rational case of Frobenius' identity \cite{F}.
Also, if we take $k=l=1$ in (\ref{schur1}), we obtain
 the rational case of an elliptic generalization of (\ref{schur})
 given in \cite{O4}.

This paper is organized as follows.
Sections 2 and 3 are devoted to the proof of Theorem~\ref{thm:main}.
In Section 2, we prove the identity (\ref{main2})
 by using the Pfaffian version of Desnanot--Jacobi formula
 and induction.
In Section~3, we give a homogeneous version of the identity (\ref{main2})
 and derive the other three identities (\ref{main1}), (\ref{main3})
 and (\ref{main4}).
A variation of the main identities is given in Section~4,
 and another instance of a Cauchy-type determinant identity is presented
 in Section~5.
Also we present a formula expressing the determinant of $V^{n,n}$
 in terms of the hyperpfaffian.
In the last section, we give an application of the identity (\ref{main2})
 to the Littlewood--Richardson coefficients.

Here we recall the definition of Pfaffians.
Given a $2n \times 2n$ skew-symmetric matrix $A=(a_{ij})_{1\leq i,j\leq2n}$,
 the Pfaffian of $A$ is defined by
$$
\Pf (A) = \sum_\sigma
 \sgn(\sigma) a_{\sigma(1),\sigma(2)} a_{\sigma(3),\sigma(4)}
 \cdots a_{\sigma(2n-1),\sigma(2n)},
$$
where $\sigma$ runs over all permutations of $2n$ letters $1, 2, \cdots, 2n$
 satisfying
\begin{gather*}
\sigma(1) < \sigma(2), \quad \sigma(3) < \sigma(4),
 \quad \cdots, \quad \sigma(2n-1) < \sigma(2n),
\\
\sigma(1) < \sigma(3) < \cdots < \sigma(2n-1).
\end{gather*}

\section{
Proof of the identity (\ref{main2}) in Theorem~\ref{thm:main}
}


In this section, we give a proof of the identity (\ref{main2}).
First we show (\ref{main2}) in the special case
 where $n=2$ by using induction.
Then we apply the Desnanot--Jacobi formula for Pfaffians
 to reduce the proof of the general case to this special case.

First we prove the case of $n=2$ by induction on $p+q+r+s$.

\begin{prop} \label{prop:n=2}
Let $p$, $q$, $r$ and $s$ be nonnegative integers.
For vectors of variables
 $\vectx$, $\vecta$, $\vectb$ of length $4$,
 $\vectz$, $\vectc$ of length $p+q$,
 and $\vectw$, $\vectd$ of length $r+s$,
 we have
\begin{multline}
\Pf \left(
 \frac{ \det V^{p+1,q+1}(x_i,x_j,\vectz ; a_i,a_j,\vectc)
        \det V^{r+1,s+1}(x_i,x_j,\vectw ; b_i,b_j,\vectd) }
      { x_j - x_i }
\right)_{1 \le i, j \le 4}
\\
=
\frac{1}{\prod_{1 \le i < j \le 4}(x_j - x_i)}
 \det V^{p,q}(\vectz ; \vectc)
 \det V^{r,s}(\vectw ; \vectd)
\\
\times
 \det V^{p+2,q+2}(\vectx,\vectz ; \vecta,\vectc)
 \det V^{r+2,s+2}(\vectx,\vectw ; \vectb,\vectd).
\label{eq:n=2}
\end{multline}
\end{prop}

In the induction step of the proof, we need relations between
 $\det V^{p,q}$ and $\det V^{p-1,q}$ (or $\det V^{q,p}$).


\begin{lemma} \label{lem:rel-v}
\begin{roster}{(2)}
\item[(1)]
If $p \ge q$ and $p \ge 1$, then we have
\begin{equation}
\det V^{p,q}(\vectx ; \vecta)
 =
\prod_{i=1}^{p+q-1}(x_{p+q}-x_i)
\cdot
\det V^{p-1,q}(x_1, \cdots, x_{p+q-1} ; a'_1, \cdots, a'_{p+q-1}),
\label{rel-v1}
\end{equation}
where we put
$$
a'_i = \frac{a_i-a_{p+q}}{x_i-x_{p+q}}
\quad(1 \le i \le p+q-1).
$$
\item[(2)]
For nonnegative integers $p$ and $q$, we have
\begin{equation}
\det V^{p,q}(\vectx ; \vecta)
=
(-1)^{pq} \prod_{i=1}^{p+q} a_i
\cdot
\det V^{q,p}(\vectx ; \vecta^{-1}),
\label{rel-v2}
\end{equation}
where $\vecta^{-1} = (a_1^{-1}, \cdots, a_{p+q}^{-1})$.
\end{roster}
\end{lemma}

\begin{demo}{Proof}
(1)
We put $m=p+q$.
By subtracting the $i$th column multiplied by $a_m$ from the $(p+i)$th column
 for $i=1, \cdots, q$,
 and
by subtracting the $i$th column multiplied by $x_m$ from the $(i+1)$th column
 for $i=p-1, \cdots, 1$,
we obtain
\begin{align*}
&
\det V^{p,q}(\vectx ; \vecta)
\\
&=
\det \begin{pmatrix}
1 & x_1-x_m & \cdots & (x_1-x_m) x_1^{p-2}
 & a_1-a_m & \cdots & (a_1-a_m) x_1^{q-1} \\
\vdots & \vdots & & \vdots
 & \vdots & & \vdots \\
1 & x_{m-1}-x_m & \cdots & (x_{m-1}-x_m) x_{m-1}^{p-2}
 & a_{m-1}-a_m & \cdots & (a_{m-1}-a_m) x_{m-1}^{q-1}
 \\
1 & 0 & \cdots & 0
 & 0 & \cdots & 0 \\
\end{pmatrix}
\\
&=
(-1)^{m+1} \prod_{k=1}^{m-1} (x_k - x_m)
\cdot
\det \begin{pmatrix}
1 & x_1 & \cdots & x_1^{p-2}
 & a'_1 & a'_1 x_1 & \cdots & a'_1 x_1^{q-1} \\
\vdots & \vdots & & \vdots
 & \vdots & \vdots & & \vdots \\
1 & x_{m-1} & \cdots & x_{m-1}^{p-2}
 & a'_{m-1} & a'_{m-1} x_{m-1} & \cdots & a'_{m-1} x_{m-1}^{q-1}
\end{pmatrix}
\\
&=
\prod_{k=1}^{p+q-1}(x_{p+q}-x_k)
\cdot
\det V^{p-1,q}(x_1, \cdots, x_{p+q-1} ; a'_1, \cdots, a'_{p+q-1}).
\end{align*}

(2)
As before we put $m=p+q$.
By performing an appropriate permutation of the columns and
 by dividing the $i$th row by $a_i$,
 we obtain
\begin{align*}
&
\det V^{p,q}(\vectx;\vecta)
\\
&\quad=
(-1)^{pq}
\det \begin{pmatrix}
a_1 & a_1 x_1 & \cdots & a_1 x_1^{q-1}
 & 1 & x_1 & \cdots & x_1^{p-1} \\
\vdots & \vdots & & \vdots
 & \vdots & \vdots & & \vdots \\
a_m & a_m x_m & \cdots & a_m x_m^{q-1}
 & 1 & x_m & \cdots & x_m^{p-1}
\end{pmatrix}
\\
&\quad=
(-1)^{pq} \prod_{k=1}^{p+q} a_k \cdot
\det \begin{pmatrix}
1 & x_1 & \cdots & x_1^{q-1}
 & a_1^{-1} & a_1^{-1} x_1 & \cdots & a_1^{-1} x_1^{p-1} \\
\vdots & \vdots & & \vdots
 & \vdots & \vdots & & \vdots \\
1 & x_m & \cdots & x_m^{q-1}
 & a_m^{-1} & a_m^{-1} x_m & \cdots & a_m^{-1} x_m^{q-1}
\end{pmatrix}
\\
&\quad=
(-1)^{pq} \prod_{k=1}^{p+q} a_k \cdot \det V^{q,p}(\vectx;\vecta^{-1}).
\end{align*}
This completes the proof of the lemma.
\end{demo}

\begin{demo}{Proof of Proposition~\ref{prop:n=2}}
We prove ({\ref{eq:n=2}}) by induction on $p+q+r+s$.
If $p+q+r+s = 0$, i.e., $p=q=r=s=0$, then one can easily check
 the equality in (\ref{eq:n=2}) by a direct computation.

Suppose $p+q+r+s > 0$.
By symmetry, we may assume $p+q > 0$ without loss of generality.
First we consider the case where $p \ge q$.
Using the relation (\ref{rel-v1}) in Lemma~\ref{lem:rel-v}, we have
\begin{multline*}
\det V^{p+1,q+1}(x_i,x_j,\vectz ; a_i,a_j,\vectc)
\\
=
(z_{p+q} - x_i)(z_{p+q} - x_j)
\prod_{k=1}^{p+q-1} (z_{p+q} - z_k)
\cdot \det V^{p,q+1}(x_i,x_j,\tilde\vectz ; a'_i,a'_j,\vectcp),
\end{multline*}
where $\tilde\vectz = (z_1, \cdots, z_{p+q-1})$ and
$a'_i$, $a'_j$ and $\vectcp = (c'_1, \cdots, c'_{p+q-1})$ are
 given by
$$
a'_k = \frac{a_k - c_{p+q}}{x_k - z_{p+q}}
\quad(k = i, j),
\quad
c'_l = \frac{c_l - c_{p+q}}{z_l - z_{p+q}}
\quad(1 \le l \le p+q-1).
$$
Hence we have
\begin{align*}
&
\Pf \left(
 \frac{ \det V^{p+1,q+1}(x_i,x_j,\vectz ; a_i,a_j,\vectc)
        \det V^{r+1,s+1}(x_i,x_j,\vectw ; b_i,b_j,\vectd) }
      { x_j - x_i }
\right)_{1 \le i, j \le 4}
\\
&\quad=
\prod_{i=1}^4 (z_{p+q} - x_i)
\prod_{i=1}^{p+q-1} (z_{p+q} - z_i)^2
\\
&\quad\quad\times
\Pf \left(
 \frac{ \det V^{p,q+1}(x_i,x_j,\tilde\vectz ; a'_i,a'_j,\vectcp)
        \det V^{r+1,s+1}(x_i,x_j,\vectw ; b_i,b_j,\vectd) }
      { x_j - x_i }
\right)_{1 \le i, j \le 4}
\intertext{by the induction hypothesis,
}
&\quad=
\prod_{i=1}^4 (z_{p+q} - x_i)
\prod_{i=1}^{p+q-1} (z_{p+q} - z_i)^2
\\
&\quad\quad\times
\frac{1}{\prod_{1 \le i < j \le 4}(x_j - x_i)}
 \det V^{p-1,q}(\tilde\vectz ; \vectcp)
 \det V^{r,s}(\vectw ; \vectd)
\\
&\quad\quad\times
 \det V^{p+1,q+2}(x_1, x_2, x_3, x_4,\tilde\vectz ;
  a'_1,a'_2,a'_3,a'_4,\vectcp)
 \det V^{r+2,s+2}(\vectx,\vectw ; \vectb,\vectd)
\\
\intertext{
by using the relation (\ref{rel-v1}) again,
}
&\quad=
\frac{1}{\prod_{1 \le i < j \le 4}(x_j - x_i)}
 \det V^{p,q}(\vectz ; \vectc)
 \det V^{r,s}(\vectw ; \vectd)
 \det V^{p+2,q+2}(\vectx,\vectz ; \vecta,\vectc)
 \det V^{r+2,s+2}(\vectx,\vectw ; \vectb,\vectd).
\end{align*}

If $p < q$, then we use the relation (\ref{rel-v2}) in Lemma~\ref{lem:rel-v}
 and the case we have just proven.
Then we see that
\begin{align*}
&
\Pf \left(
 \frac{ \det V^{p+1,q+1}(x_i,x_j,\vectz ; a_i,a_j,\vectc)
        \det V^{r+1,s+1}(x_i,x_j,\vectw ; b_i,b_j,\vectd) }
      { x_j - x_i }
\right)_{1 \le i, j \le 4}
\\
&\quad=
\Pf \left(
 \frac{ a_i a_j \prod_{k=1}^{p+q} c_k\cdot
        \det V^{q+1,p+1}(x_i,x_j,\vectz ; a_i^{-1},a_j^{-1},\vectc^{-1})
        \det V^{r+1,s+1}(x_i,x_j,\vectw ; b_i,b_j,\vectd) }
      { x_j - x_i }
\right)_{1 \le i, j \le 4}
\\
&\quad=
\prod_{i=1}^4 a_i \prod_{k=1}^{p+q} c_k^2
\\
&\quad\quad\times
\frac{1}{\prod_{1 \le i < j \le 2n}(x_j - x_i)}
 \det V^{q,p}(\vectz ; \vectc^{-1})
 \det V^{r,s}(\vectw ; \vectd)
 \det V^{q+2,p+2}(\vectx,\vectz ; \vecta^{-1},\vectc^{-1})
 \det V^{r+2,s+2}(\vectx,\vectw ; \vectb,\vectd)
\\
&\quad=
\frac{1}{\prod_{1 \le i < j \le 4}(x_j - x_i)}
 \det V^{p,q}(\vectz ; \vectc)
 \det V^{r,s}(\vectw ; \vectd)
 \det V^{p+2,q+2}(\vectx,\vectz ; \vecta,\vectc)
 \det V^{r+2,s+2}(\vectx,\vectw ; \vectb,\vectd).
\end{align*}
This completes the proof of Proposition~\ref{prop:n=2}.
\end{demo}


Here we recall the Desnanot--Jacobi formula for determinants and
 Pfaffians.
Given a square matrix $A$ and indices $i_1, \cdots, i_r, j_1, \cdots, j_r$,
 we denote by $A^{i_1, \cdots, i_r}_{j_1, \cdots, j_r}$ the matrix
 obtained by removing the rows $i_1, \cdots, i_r$ and
 the columns $j_1, \cdots, j_r$ of $A$.

\begin{lemma} \label{lem:dodgson}
\begin{roster}{(2)}
\item[(1)]
If $A$ is a square matrix, then we have
\begin{equation}
\det A^1_1 \cdot \det A^2_2 - \det A^1_2 \cdot \det A^2_1
 = \det A \cdot \det A^{1,2}_{1,2}.
\label{det-dodgson}
\end{equation}
\item[(2)]
If $A$ is a skew-symmetric matrix $A$, then we have
\begin{equation}
\Pf A^{1,2}_{1,2} \cdot \Pf A^{3,4}_{3,4}
 - \Pf A^{1,3}_{1,3} \cdot \Pf A^{2,4}_{2,4}
 + \Pf A^{1,4}_{1,4} \cdot \Pf A^{2,3}_{2,3}
 =
 \Pf A \cdot \Pf A^{1,2,3,4}_{1,2,3,4}.
\label{pf-dodgson}
\end{equation}
\end{roster}
\end{lemma}

This Pfaffian analogue of the Desnanot-Jacobi formula is
 given in \cite{Kn} and \cite{IW2}.


If $\vectx=(x_{1},\dots,x_{n})$ is a vector of variables
 and $1 \le i_1 < \dots < i_r \le n$ are indices,
we denote by $\vectx^{(i_1,\cdots,i_r)}$ the vector obtained from $\vectx$
 by removing the variables $x_{i_1}, \cdots, x_{i_r}$.

\begin{demo}{Proof of the identity in (\ref{main2}) in Theorem~\ref{thm:main}}
We proceed by induction on $n$.
If $n=1$, then there is nothing to prove,
 and, if $n = 2$, we already proved (\ref{main2}) in
 Proposition~\ref{prop:n=2}.

Suppose $n \ge 3$.
Apply the Desnanot--Jacobi formula for Pfaffians (\ref{pf-dodgson})
 in Lemma~\ref{lem:dodgson} to the skew-symmetric matrix
$$
A =
\left(
 \frac{ \det V^{p+1,q+1}(x_i,x_j,\vectz ; a_i,a_j,\vectc)
        \det V^{r+1,s+1}(x_i,x_j,\vectw ; b_i,b_j,\vectd) }
      { x_j - x_i }
\right)_{1 \le i, j \le 2n}.
$$
Then the induction hypothesis tells us that, for $1 \le k < l \le 4$,
 we have
\begin{align*}
\Pf A^{k,l}_{k,l}
&=
\frac{ 1 }
     { (x_{l'} - x_{k'}) \prod_{i=5}^{2n} (x_i - x_{k'})(x_i - x_{l'})
       \prod_{5 \le i < j \le 2n}(x_j - x_i) }
\\
&\quad\times
 \det V^{p,q}(\vectz ; \vectc)^{n-2}
 \det V^{r,s}(\vectw ; \vectd)^{n-2}
\\
&\quad\times
 \det V^{n+p-1,n+q-1}(\vectx^{(k,l)},\vectz ; \vecta^{(k,l)},\vectc)
 \det V^{n+r-1,n+s-1}(\vectx^{(k,l)},\vectw ; \vectb^{(k,l)},\vectd),
\end{align*}
where $k'$ and $l'$ are the indices satisfying $\{ k, l, k',l' \}
 = \{ 1, 2, 3, 4 \}$ and $k < l$, $k' < l'$,
and
\begin{align*}
&
\Pf A^{1,2,3,4}_{1,2,3,4}
\\
&\quad=
\frac{ 1 }
     { \prod_{5 \le i < j \le 2n}(x_j - x_i) }
 \det V^{p,q}(\vectz ; \vectc)^{n-3}
 \det V^{r,s}(\vectw ; \vectd)^{n-3}
\\
&\quad\quad\times
 \det V^{n+p-2,n+q-2}(\vectx^{(1,2,3,4)},\vectz ; \vecta^{(1,2,3,4)},\vectc)
 \det V^{n+r-2,n+s-2}(\vectx^{(1,2,3,4)},\vectw ; \vectb^{(1,2,3,4)},\vectd).
\end{align*}
Hence, by applying (\ref{pf-dodgson}) and cancelling the common factors,
 we see that, in order to prove (\ref{main2}), it suffices to show
\begin{align*}
&\quad
 \frac{ \det V^{n+p-1,n+q-1}(\vectx^{(1,2)},\vectz ; \vecta^{(1,2)},\vectc)
        \det V^{n+r-1,n+s-1}(\vectx^{(1,2)},\vectw ; \vectb^{(1,2)},\vectd) }
      { x_2 - x_1 }
\\
&\quad\times
 \frac{ \det V^{n+p-1,n+q-1}(\vectx^{(3,4)},\vectz ; \vecta^{(3,4)},\vectc)
        \det V^{n+r-1,n+s-1}(\vectx^{(3,4)},\vectw ; \vectb^{(3,4)},\vectd) }
      { x_4 - x_3 }
\\
&\quad-
 \frac{ \det V^{n+p-1,n+q-1}(\vectx^{(1,3)},\vectz ; \vecta^{(1,3)},\vectc)
        \det V^{n+r-1,n+s-1}(\vectx^{(1,3)},\vectw ; \vectb^{(1,3)},\vectd) }
      { x_3 - x_1 }
\\
&\quad\times
 \frac{ \det V^{n+p-1,n+q-1}(\vectx^{(2,4)},\vectz ; \vecta^{(2,4)},\vectc)
        \det V^{n+r-1,n+s-1}(\vectx^{(2,4)},\vectw ; \vectb^{(2,4)},\vectd) }
      { x_4 - x_2 }
\\
&\quad+
 \frac{ \det V^{n+p-1,n+q-1}(\vectx^{(1,4)},\vectz ; \vecta^{(1,4)},\vectc)
        \det V^{n+r-1,n+s-1}(\vectx^{(1,4)},\vectw ; \vectb^{(1,4)},\vectd) }
      { x_4 - x_1 }
\\
&\quad\times
 \frac{ \det V^{n+p-1,n+q-1}(\vectx^{(2,3)},\vectz ; \vecta^{(2,3)},\vectc)
        \det V^{n+r-1,n+s-1}(\vectx^{(2,3)},\vectw ; \vectb^{(2,3)},\vectd) }
      { x_3 - x_2 }
\\
&
=
 \det V^{n+p-2,n+q-2}(\vectx^{(1,2,3,4)},\vectz ; \vecta^{(1,2,3,4)},\vectc)
 \det V^{n+r-2,n+s-2}(\vectx^{(1,2,3,4)},\vectw ; \vectb^{(1,2,3,4)},\vectd)
\\
&\quad\times
 \frac{ \det V^{n+p,n+q}(\vectx,\vectz ; \vecta,\vectc)
        \det V^{n+r,n+s}(\vectx,\vectw ; \vectb,\vectd) }
      { \prod_{1 \le i < j \le 4} ( x_j - x_i ) }.
\end{align*}
This is equivalent to the identity (\ref{eq:n=2}) with
 $\vectz$, $\vectc$, $\vectw$, $\vectd$ replaced by
$$
\vectz \leftarrow (\vectx^{(1,2,3,4)}, \vectz),
\quad
\vectc \leftarrow (\vecta^{(1,2,3,4)}, \vectc),
\quad
\vectw \leftarrow (\vectx^{(1,2,3,4)}, \vectw),
\quad
\vectd \leftarrow (\vectb^{(1,2,3,4)}, \vectd),
$$
respectively, and it is proven in Proposition~\ref{prop:n=2}.
This completes the proof of (\ref{main2}).
\end{demo}

\begin{remark}
We can also reduce the proof of the other identities (\ref{main1}),
 (\ref{main3}) and (\ref{main4}) in Theorem~\ref{thm:main}
 to the case of $n=2$ with the help of the Desnanot-Jacobi formulae.
It is easy to show the case of $n=2$ of (\ref{main1}) by using
 the relations in Lemma~\ref{lem:rel-v} and induction on $p+q$.
Also we can prove (\ref{main3}) (resp. (\ref{main4})) in the case of $n=2$,
 by regarding both sides as polynomials in $z_p$
 and showing that the values coincide at $(2p+3)$ distinct points
 $z_1, \cdots, z_{p-1}, z_1^{-1}, \cdots, z_{p-1}^{-1} x_1, x_2, y_1, y_2$
 and $-1$
 (resp. $z_1, \cdots, z_{p-1}, z_1^{-1}, \cdots, z_{p-1}^{-1} x_1, x_2,
 x_3, x_4$ and $-1$) with the help of induction.

But, in this paper, we adopt another method, namely,
we \lq\lq homogenize\rq\rq\  the identity (\ref{main2}) and
 derive the other identities from this homogeneous version (\ref{homog1}).
\end{remark}

\section{
Proof of the identities (\ref{main1}), (\ref{main2}) and (\ref{main3})
}

In this section, we give a homogeneous version of the identity (\ref{main2}),
 which is shown in the previous section,
 and derive the identities (\ref{main1}), (\ref{main3}) and (\ref{main4})
 from this homogeneous version.

Throughout the remaining of this paper, we use the following notation for
 vectors $\vectx = (x_1, \cdots, x_n)$ and $\vecty = (y_1, \cdots, y_n)$:
$$
\vectx + \vecty
 = (x_1 + y_1, \ldots, x_n + y_n),
\quad
\vectx \vecty
 = (x_1 y_1, \ldots, x_n y_n),
$$
and, for integers $k$ and $l$,
$$
\vectx^k
 = (x_1^k, \ldots, x_n^k),
\quad
\vectx^k \vecty^l
 = (x_1^k y_1^l, \ldots, x_n^k y_n^l).
$$

We introduce a homogeneous version of the matrix $V^{p,q}(\vectx;\vecta)$
 as follows.
For vectors $\vectx$, $\vecty$, $\vecta$, $\vectb$ of length $n$
 and nonnegative integers $p$, $q$ with $p+q = n$, we define
 a matrix $\U{p}{q}{\vectx}{\vecty}{\vecta}{\vectb}$ to be
 the $n \times n$ matrix with $i$th row
$$
(a_i x_i^{p-1}, a_i x_i^{p-2} y_i, \cdots, a_i y_i^{p-1},
 b_i x_i^{q-1}, b_i x_i^{q-2} y_i, \cdots, b_i y_i^{q-1}).
$$
Then the following relations between $\det V^{p,q}$ and $\det U^{p,q}$
 are easily shown by elementary transformations, so we omit their proofs.

\begin{lemma} \label{lem:rel-uv}
\begin{equation}
\det\U{p}{q}{\vectx}{\vecty}{\vecta}{\vectb}
 =
\prod_{k=1}^{p+q} a_k x_k^{p-1}
 \cdot 
\det V^{p,q} \left( \vectx^{-1} \vecty ; \vecta^{-1} \vectb \vectx^{q-p} \right).
\label{rel-uv1}
\end{equation}
In particular,
\begin{equation}
\det V^{p,q}(\vectx;\vecta)
 =
\det\U{p}{q}{\vectone}{\vectx}{\vectone}{\vecta},
\label{rel-uv2}
\end{equation}
where $\vectone=(1,\dots,1)$.
\end{lemma}

Now we give a homogeneous version of the identity (\ref{main2}).

\begin{theorem} \label{thm:general}
Let $n$ be a positive integer
and let $p$, $q$, $r$ and $s$ be nonnegative integers.
Suppose that
the vectors $\vectx$, $\vecty$, $\vecta$, $\vectb$, $\vectc$, $\vectd$
 have length $2n$,
the vectors $\vectxi$, $\vecteta$, $\vectalpha$, $\vectbeta$ have
 length $p+q$,
and the vectors $\vectzeta$, $\vectomega$, $\vectgamma$, $\vectdelta$
 have length $r+s$.
Then we have
\begin{multline}
\Pf \left(
 \frac{ \det \U{p+1}{q+1}{x_i,x_j,\vectxi}{y_i,y_j,\vecteta}
               {a_i,a_j,\vectalpha}{b_i,b_j,\vectbeta}
        \det \U{r+1}{s+1}{x_i,x_j,\vectzeta}{y_i,y_j,\vectomega}
               {c_i,c_j,\vectgamma}{d_i,d_j,\vectdelta}
      }
      { \det \begin{pmatrix} x_i & x_j \\ y_i & y_j \end{pmatrix}
      }
\right)_{1 \le i < j \le 2n}
\\
=
\frac{ 1 }
     { \prod_{1 \le i < j \le 2n}
        \det \begin{pmatrix} x_i & x_j \\ y_i & y_j \end{pmatrix}
     }
\det \U{p}{q}{\vectxi}{\vecteta}{\vectalpha}{\vectbeta}^{n-1}
\det \U{r}{s}{\vectzeta}{\vectomega}{\vectgamma}{\vectdelta}^{n-1}
\\
\times
\det \U{n+p}{n+q}{\vectx,\vectxi}{\vecty,\vecteta}
       {\vecta,\vectalpha}{\vectb,\vectbeta}
\det \U{n+r}{n+s}{\vectx,\vectzeta}{\vecty,\vectomega}
       {\vectc,\vectgamma}{\vectd,\vectdelta}.
\label{homog1}
\end{multline}
\end{theorem}

The special case of $p=q=r=s=0$ of this identity (\ref{homog1})
 is given by M.~Ishikawa \cite[Theorem~3.1]{I},
 and is one of the key ingredients of his proof of Stanley's conjecture.

The identity (\ref{main2}) is equivalent to (\ref{homog1}).
It follows from the relation (\ref{rel-uv2}) that (\ref{homog1}) specialize
 to (\ref{main2}) by setting
$$
\vectx = \vecta = \vectc = \vectone_{2n},
\quad
\vectxi = \vectalpha = \vectone_{p+q},
\quad
\vectzeta = \vectgamma = \vectone_{r+s}
$$
and renaming the variables.
As we see in the following proof, we derive (\ref{homog1}) from (\ref{main2}).

\begin{demo}{Proof of Theorem~\ref{thm:general}}
In the identity (\ref{main2}), we substitute as follows:
\begin{equation}
\begin{gathered}
\vectx \leftarrow \vectx^{-1} \vecty,
\quad
\vectz \leftarrow \vectxi^{-1} \vecteta,
\quad
\vectw \leftarrow \vectzeta^{-1} \vectomega,
\\
\vecta \leftarrow \vecta^{-1} \vectb \vectx^{q-p},
\quad
\vectc \leftarrow \vectalpha^{-1} \vectbeta \vectxi^{q-p},
\quad
\vectb \leftarrow \vectc^{-1} \vectd \vectx^{s-r},
\quad
\vectd \leftarrow \vectgamma^{-1} \vectdelta \vectzeta^{s-r}.
\end{gathered}
\label{subs1}
\end{equation}
By using the relations (\ref{rel-uv1}) and 
$$
\frac{y_j}{x_j} - \frac{y_i}{x_i}
 =
x_i^{-1} x_j^{-1}
\det \begin{pmatrix} x_i & x_j \\ y_i & y_j \end{pmatrix},
$$
we see that the $(i,j)$ entry of the left-hand side of (\ref{main2})
 under the substitution (\ref{subs1}) becomes
\begin{multline*}
\frac{
 \left(
  a_i a_j x_i^{p} x_j^{p}
  \prod_{k=1}^{p+q} \alpha_k \xi_k^{p}
 \right)^{-1}
 \left(
  c_i c_j x_i^{r} x_j^{r}
  \prod_{k=1}^{r+s} \gamma_k \zeta_k^{r}
 \right)^{-1}
}{
 x_i^{-1} x_j^{-1}
}
\\
\times
\frac{
 \det \U{p+1}{q+1}{x_i,x_j,\vectxi}{y_i,y_j,\vecteta}
        {a_i,a_j,\vectalpha}{b_i,b_j,\vectbeta}
 \det \U{r+1}{s+1}{x_i,x_j,\vectzeta}{y_i,y_j,\vectomega}
        {c_i,c_j,\vectgamma}{d_i,d_j,\vectdelta}
}{
 \det \begin{pmatrix} x_i & x_j \\ y_i & y_j \end{pmatrix}
}.
\end{multline*}
Hence, by noting the linearity of Pfaffian
$$
\Pf \left( \lambda \alpha_i \alpha_j a_{ij} \right)_{1 \le i < j \le 2n}
 = 
\lambda^n \prod_{i=1}^{2n} \alpha_i \cdot
 \Pf \left( a_{ij} \right)_{1 \le i < j \le 2n},
$$
we can see the Pfaffian on the left-hand side of (\ref{main2}) becomes
\begin{multline*}
\left(
 \prod_{i=1}^{2n} a_i c_i x_i^{p+r-1}
\right)^{-1}
\left(
 \prod_{i=1}^{p+q} \alpha_i \xi_i^p
 \prod_{i=1}^{r+s} \gamma_i \zeta_i^r
\right)^{-n}
\\
\times
\Pf \left(
\frac{
 \det \U{p+1}{q+1}{x_i,x_j,\vectxi}{y_i,y_j,\vecteta}
        {a_i,a_j,\vectalpha}{b_i,b_j,\vectbeta}
 \det \U{r+1}{s+1}{x_i,x_j,\vectzeta}{y_i,y_j,\vectomega}
        {c_i,c_j,\vectgamma}{d_i,d_j,\vectdelta}
}{
 \det \begin{pmatrix} x_i & x_j \\ y_i & y_j \end{pmatrix}
}
\right)_{1 \le i < j \le 2n}.
\end{multline*}

On the other hand, using the relations (\ref{rel-uv1}) and
$$
\prod_{1 \le i < j \le 2n}
 \left( \frac{y_j}{x_j} - \frac{y_i}{x_i} \right)
=
\prod_{k=1}^{2n} x_k^{-2n+1}
\prod_{1 \le i < j \le 2n}
 \det \begin{pmatrix} x_i & x_j \\ y_i & y_j \end{pmatrix},
$$
we see that the right-hand side of (\ref{main2}) becomes
\begin{align*}
&
\frac{ 1 }
     { \left( \prod_{k=1}^{2n} x_k \right)^{-2n+1} }
\cdot
\frac{ 1 }
     { \prod_{1 \le i < j \le 2n}
         \det \begin{pmatrix} x_i & x_j \\ y_i & y_j \end{pmatrix}
     }
\\
&\quad\times
 \left( \prod_{i=1}^{p+q} \alpha_i \xi_i^{p-1} \right)^{-n+1}
 \left( \prod_{i=1}^{r+s} \gamma_i \zeta_i^{r-1} \right)^{-n+1}
 \det \U{p}{q}{\vectxi}{\vecteta}{\vectalpha}{\vectbeta}^{n-1}
 \det \U{r}{s}{\vectzeta}{\vectomega}{\vectgamma}{\vectdelta}^{n-1}
\\
&\quad\times
 \left(
   \prod_{i=1}^{2n} a_i x_i^{n+p-1}
   \prod_{i=1}^{p+q} \alpha_i \xi_i^{n+p-1}
 \right)^{-1}
 \left(
   \prod_{i=1}^{2n} c_i x_i^{n+r-1}
   \prod_{i=1}^{r+s} \gamma_i \zeta_i^{n+r-1}
 \right)^{-1}
\\
&\quad\times
 \det \U{n+p}{n+q}{\vectx,\vectxi}{\vecty,\vecteta}
        {\vecta,\vectalpha}{\vectb,\vectbeta}
 \det \U{n+r}{n+s}{\vectx,\vectzeta}{\vecty,\vectomega}
        {\vectc,\vectgamma}{\vectd,\vectdelta}.
\end{align*}
Comparing the both sides and cancelling the common factors,
 we obtain the desired identity (\ref{homog1}).
\end{demo}

In this setting, a homogeneous version of (\ref{main1}) is
 a direct consequence of (\ref{homog1}).
A key is the following relation between determinant and Pfaffian.
If $A$ is any $m \times (2n-m)$ matrix, then we have
\begin{equation}
\Pf \begin{pmatrix}
 O & A \\ - \trans A & O
\end{pmatrix}
 = \begin{cases}
(-1)^{n(n-1)/2} \det A &\text{if $m=n$}, \\
0 &\text{if $m \neq n$}.
\end{cases}
\label{pf-det}
\end{equation}

\begin{corollary} \label{cor:general}
Let $n$ be a positive integer
and let $p$ and $q$ be fixed nonnegative integers.
For vectors $\vectx$, $\vecty$, $\vectz$, $\vectw$, $\vecta$, $\vectb$,
 $\vectc$, $\vectd$ of length $n$,
 and vectors $\vectxi$, $\vecteta$, $\vectalpha$, $\vectbeta$ of length $p+q$,
 we have
\begin{multline}
\det \left(
 \frac{ \det \U{p+1}{q+1}{x_i,z_j,\vectxi}{y_i,w_j,\vecteta}
               {a_i,c_j,\vectalpha}{b_i,d_j,\vectbeta}
      }
      { \det \begin{pmatrix} x_i & z_j \\ y_i & w_j \end{pmatrix}
      }
\right)_{1 \le i, j \le n}
\\
=
\frac{ (-1)^{n(n-1)/2} }
     { \prod_{1 \le i, j \le n}
        \det \begin{pmatrix} x_i & z_j \\ y_i & w_j \end{pmatrix}
     }
\det \U{p}{q}{\vectxi}{\vecteta}{\vectalpha}{\vectbeta}^{n-1}
\det \U{n+p}{n+q}{\vectx,\vectz,\vectxi}{\vecty,\vectw,\vecteta}
       {\vecta,\vectc,\vectalpha}{\vectb,\vectd,\vectbeta}.
\label{homog2}
\end{multline}
\end{corollary}

\begin{demo}{Proof}
In (\ref{homog1}), we take $r = s = 0$ and put
\begin{equation}
\begin{gathered}
c_1 = \dots = c_n = 1, \qquad c_{n+1} = \dots = c_{2n} = 0,
\\
d_1 = \dots = d_n = 0, \qquad d_{n+1} = \dots = d_{2n} = 1.
\end{gathered}
\label{subs2}
\end{equation}
Under this substitution (\ref{subs2}), we have
$$
\det \U{1}{1}{x_i,x_j}{y_i,y_j}{c_i,c_j}{d_i,d_j}
=
\det \begin{pmatrix} c_i & d_i \\ c_j & d_j \end{pmatrix}
=
\begin{cases}
 0 &\text{if $1 \le i, j \le n$ or $n+1 \le i, j \le 2n$,}
\\
 1 &\text{if $1 \le i \le n$ and $n+1 \le j \le 2n$,} \\
-1 &\text{if $n+1 \le i \le 2n$ and $1 \le j \le n$.}
\end{cases}
$$
Hence, by (\ref{pf-det}), we see that the left-hand side of
 (\ref{homog1}) becomes
$$
(-1)^{n(n-1)/2}
\det \left(
 \frac{ \det \U{p+1}{q+1}{x_{i},x_{n+j},\vectxi}{y_{i},y_{n+j},\vecteta}
         {a_i,a_{n+j},\vectalpha}{b_i,b_{n+j},\vectbeta}
      }
      { \det \begin{pmatrix} x_i & x_{n+j} \\ y_i & y_{n+j} \end{pmatrix}
      }
\right)_{1 \le i,j \le n}.
$$
On the other hand, under the specialization (\ref{subs2}), we have
\begin{align*}
\det \U{n}{n}{\vectx}{\vecty}{\vectc}{\vectd}
 &=
\det \begin{pmatrix}
 \left( x_i^{n-j} y_i^{j-1} \right)_{1 \le i, j \le n} & O \\
 O & \left( x_{n+i}^{n-j} y_{n+i}^{j-1} \right)_{1 \le i, j \le n}
\end{pmatrix}
\\
 &=
\prod_{1 \le i < j \le n}
 \det \begin{pmatrix} x_i & x_j \\ y_i & y_j \end{pmatrix}
\prod_{n+1 \le i < j \le 2n}
 \det \begin{pmatrix} x_i & x_j \\ y_i & y_j \end{pmatrix}.
\end{align*}
Thus the right-hand side of (\ref{homog1}) becomes
$$
\frac{ 1 }
     { \prod_{1 \le i, j \le n}
        \det \begin{pmatrix} x_i & x_{n+j} \\ y_i & y_{n+j} \end{pmatrix}
     }
\det \U{p}{q}{\vectxi}{\vecteta}{\vectalpha}{\vectbeta}^{n-1}
\det \U{n+p}{n+q}{\vectx,\vectxi}{\vecty,\vecteta}
       {\vecta,\vectalpha}{\vectb,\vectbeta}.
$$
Lastly, if we replace the variables as $x_{n+i}=z_{i}$,
 $y_{n+i}=w_{i}$, $a_{n+i}=c_{i}$ and $b_{n+i}=d_{i}$
 for $1 \le i \le n$,
 then we obtain the desired identity (\ref{homog2}).
This completes our proof.
\end{demo}

Now it is easy to derive (\ref{main1}) from (\ref{pf-det})
 by using the relation (\ref{rel-uv1}).
To prove the remaining identities (\ref{main3}) and (\ref{main4})
 in Theorem~\ref{thm:main}, we need the following lemma.

\begin{lemma} \label{lem:rel-uw}
Let $n$ be a nonnegative integer.
\begin{roster}{(2)}
\item[(1)]
For vectors $\vectx$, $\vecta$ and $\vectone=(1, \ldots, 1)$ of length $2n$,
 we have
\begin{equation}
\det \U{n}{n}{\vectx}{\vectone+\vectx^2}{\vectone+\vecta\vectx}{\vectx+\vecta}
=
(-1)^{n(n-1)/2} \det W^{2n}(\vectx;\vecta).
\label{rel-uw1}
\end{equation}
\item[(2)]
For vectors $\vectx$, $\vecta$ and $\vectone=(1, \ldots, 1)$ of length $2n+1$,
 we have
\begin{equation}
\det \U{n}{n+1}{\vectx}{\vectone+\vectx^2}
       {\vectone+\vecta\vectx^2}{\vectone+\vecta}
=
(-1)^{n(n-1)/2} \det W^{2n+1}(\vectx;\vecta).
\label{rel-uw2}
\end{equation}
\end{roster}
\end{lemma}

\begin{demo}{Proof}
By definition, we have
\begin{multline*}
\det \U{p}{q}{\vectx}{\vectone+\vectx^2}{\vectc}{\vectd}
\\
=
\det \left(
\begin{cases}
 c_i x_i^{p-j} (1+x_i^2)^{j-1}
  &\text{if $1 \le i \le p+q$ and $1 \le j \le p$,}\\
 d_i x_i^{p+q-j} (1+x_i^2)^{j-p-1}
  &\text{if $1 \le i \le p+q$ and $p+1 \le j \le p+q$.}
\end{cases}
\right).
\end{multline*}
By performing appropriate elementary column transformations,
one can see that this determinant is equal to
\begin{equation}
\det \left(
\begin{cases}
c_i x_i^{p-1}
 &\text{if $1 \le i \le p+q$ and $j = 1$,} \\
c_i x_i^{p-j} \left( 1+x_i^{2(j-1)} \right)
 &\text{if $1 \le i \le p+q$ and $2 \le j \le p$,} \\
d_i x_i^{q-1}
 &\text{if $1 \le i \le p+q$ and $j = p+1$,} \\
d_i x_i^{p+q-j} \left( 1+x_i^{2(j-p-1)} \right)
 &\text{if $1 \le i \le p+q$ and $p+2 \le j \le p+q$.}
\end{cases}
\right).
\label{det-u}
\end{equation}

(1)
First put $q = p = n$, $c_i = 1 + a_i x_i$ and $d_i = x_i + a_i$
 for $1 \le i \le 2n$.
Then the above determinant (\ref{det-u}) is equal to
$$
\det \left(
\begin{cases}
x_i^{n-1} + a_i x_i^n
 &\text{if $1 \le i \le 2n$ and $j = 1$,} \\
x_i^{n-j} + a_i x_i^{n+j-1} + x_i^{n+j-2} + a_i x_i^{n-j+1}
 &\text{if $1 \le i \le 2n$ and $2 \le j \le n$,} \\
x_i^n + a_i x_i^{n-1}
 &\text{if $1 \le i \le 2n$ and $j = n+1$,} \\
x_i^{2n-j+1} + a_i x_i^{j-2} + x_i^{j-1} + a_i x_i^{2n-j}
 &\text{if $1 \le i \le 2n$ and $n+2 \le j \le 2n$.}
\end{cases}
\right).
$$
Subtract the first column from the $(n+2)$th column and subtract
the $(n+1)$th column from the second column, then subtract the
second column from the $(n+3)$th column and subtract the $(n+2)$th
column from the third column, and so on.
We continue these elementary column transformations until we obtain
$$
\det \left(
\begin{cases}
x_i^{n-j} + a_i x_i^{n+j-1}
 &\text{if $1 \le i \le 2n$ and $1 \le j \le n$,} \\
x_i^{j-1} + a_i x_i^{2n-j}
 &\text{if $1 \le i \le 2n$ and $n+1 \le j \le 2n$.}
\end{cases}
\right),
$$
which is $(-1)^{n(n-1)/2} \det W^{2n}(\vectx;\vecta)$.

(2)
Next we take $p=n$, $q=n+1$, $c_i = 1 + a_i x_i^2$ and
$d_i = 1 + a_i$ ($1 \le i \le 2n+1$) in (\ref{det-u}),
then we obtain
$$
\det \left(
\begin{cases}
x_i^{n-1} + a_i x_i^{n+1}
 &\text{if $1 \le i \le 2n+1$ and $j = 1$,} \\
x_i^{n-j} + a_i x_i^{n+j} + x_i^{n+j-2} + a_i x_i^{n-j+2}
 &\text{if $1 \le i \le 2n+1$ and $2 \le j \le n$,} \\
x_i^n + a_i x_i^n
 &\text{if $1 \le i \le 2n+1$ and $j = n+1$,} \\
x_i^{2n+1-j} + a_i x_i^{j-1} + x_i^{j-1} + a_i x_i^{2n+1-j}
 &\text{if $1 \le i \le 2n+1$ and $n+2 \le j \le 2n+1$.}
\end{cases}
\right).
$$
We subtract the first column from the $(n+2)$th column
and subtract the $(n+1)$th column from the second column,
and then subtract the second column from the $(n+3)$th column
and subtract the $(n+2)$th column from the third column, and so on.
We continue these elementary transformations until we obtain
$$
\det \left(
\begin{cases}
x_i^{n-i} + a_i x_i^{n+j}
 &\text{if $1 \le i \le 2n+1$ and $1 \le j \le n$,} \\
x_i^{j-1} + a_i x_i^{2n+1-j}
 &\text{if $1 \le i \le 2n+1$ and $n+1 \le j \le 2n+1$.}
\end{cases}
\right),
$$
which is equal to $(-1)^{n(n-1)/2} \det W^{2n+1}(\vectx;\vecta)$.
This completes our proof.
\end{demo}

Now we can finish our proof of Theorem~\ref{thm:main}.

\begin{demo}{Proof of the identities (\ref{main1}), (\ref{main3})
 and (\ref{main4}) in Theorem~\ref{thm:main}}
As we mentioned before, the identity (\ref{main1}) follows
 from (\ref{homog2}) by virtue of (\ref{rel-uv1}).

We derive (\ref{main4}) from (\ref{homog1}).
First we consider the case where both $p=2l$ and $q=2m$ are even.
In (\ref{homog1}), we take $p=q=l$ and $r=s=m$, and perform
 the following substitutions:
\begin{gather*}
\vectx \leftarrow \vectx,
\quad
\vecty \leftarrow \vectone + \vectx^2,
\quad
\vectxi \leftarrow \vectz,
\quad
\vecteta \leftarrow \vectone + \vectz^2,
\quad
\vectzeta \leftarrow \vectw,
\quad
\vectomega \leftarrow \vectone + \vectw^2,
\\
\vecta \leftarrow \vectone + \vecta \vectx,
\quad
\vectb \leftarrow \vectx + \vecta,
\quad
\vectc \leftarrow \vectone + \vectb \vectx,
\quad
\vectd \leftarrow \vectx + \vectb,
\\
\vectalpha \leftarrow \vectone + \vectc \vectz,
\quad
\vectbeta \leftarrow \vectz + \vectc,
\quad
\vectgamma \leftarrow \vectone + \vectd \vectw,
\quad
\vectdelta \leftarrow \vectw + \vectd.
\end{gather*}
By using the relation
$$
\det \begin{pmatrix} x_i & x_j \\ 1 + x_i^2 & 1 + x_j^2 \end{pmatrix}
=
(x_i - x_j)(1 - x_i x_j),
$$
and (\ref{rel-uw1}) in Lemma~\ref{lem:rel-uw}, we see that
 the identity (\ref{homog1}) becomes
\begin{align*}
&
\Pf \left(
\frac{ (-1)^{l(l+1)/2 + m(m+1)/2}
       \det W^{2l+2}(x_i,x_j,\vectz ; a_i,a_j,\vectc)
       \det W^{2m+2}(x_i,x_j,\vectw ; b_i,b_j,\vectd)
     }
     { (x_i - x_j)(1 - x_i x_j) }
\right)_{1 \le i < j \le 2n}
\\
&\quad=
(-1)^{(n-1) \binom{l}2 + (n-1) \binom{m}2 + \binom{n+l}2 + \binom{n+m}2}
\\
&\quad\quad\times
\frac{ 1 }
     { \prod_{1 \le i < j \le 2n}
       (x_i - x_j)(1 - x_i x_j)
     }
\\
&\quad\quad\times
\det W^{2l}(\vectz;\vectc)^{n-1}
\det W^{2m}(\vectw;\vectd)^{n-1}
\det W^{2n+2l}(\vectx,\vectz;\vecta,\vectc)
\det W^{2n+2m}(\vectx,\vectw;\vectb,\vectd).
\end{align*}
Since we have
\begin{multline*}
(n-1) l(l-1)/2 + (n-1) m(m-1)/2 + (n+l)(n+l-1)/2 + (n+m)(n+m-1)/2
\\
= n \left\{ l(l+1)/2 + m(m+1)/2 \right\} + n(n-1),
\end{multline*}
we obtain the identity (\ref{main4}) when $p$ and $q$ are both even.
We can prove the other cases similarly by using (\ref{rel-uw1})
 or (\ref{rel-uw2}) according as $p$ and $q$ are even or odd.

Also the remaining identity (\ref{main3}) can be derived from
 (\ref{homog2}) by using (\ref{rel-uw1}) or (\ref{rel-uw2}).
The details are left to the reader.
\end{demo}


\section{
A variation of the determinant and Pfaffian identities
}

In this section, we give a variation of the identities in
 Theorem~\ref{thm:main}, which can be regarded as a generalization
 of an identity of T.~Sundquist \cite{Su}.
This variation is proposed by one of the authors.

Let $n$ be a positive integer and let $p$ and $q$ be nonnegative integers
 with $p+q = n$.
Let $\vectx = (x_1, \cdots, x_n)$ and
 $\vecta = (a_1, \cdots, a_n)$ be vectors of variables.
For partitions $\lambda$ and $\mu$ with $l(\lambda) \le p$ and $l(\mu) \le q$,
 we define a matrix $V^{p,q}_{\lambda,\mu}(\vectx ; \vecta)$ to be
 the $n \times n$ matrix with $i$th row
$$
(x_i^{\lambda_p}, x_i^{\lambda_{p-1}+1}, x_i^{\lambda_{p-2}+2},
 \cdot, x_i^{\lambda_1+p-1},
a_i x_i^{\mu_q}, a_i x_i^{\mu_{q-1}+1}, a_i x_i^{\mu_{q-2}+2}, \cdots,
 a_i x_i^{\mu_1+q-1}).
$$
For example, if $\lambda = \mu = \emptyset$, then we have
 $V^{p,q}_{\emptyset,\emptyset}(\vectx;\vecta) = V^{p,q}(\vectx;\vecta)$.
Let $P_n$ denote the set of integer partitions of the form
$(\alpha_1,\ldots,\alpha_r|\alpha_1+1,\ldots, \alpha_r+1)$ in the
Frobenius notation with $\alpha_1+2\leq n$. (So $n\geq 2$).
We define
$$
F^{p,q}(\vectx ;\vecta)
 =
\sum_{\lambda \in \Par_p, \, \mu \in \Par_q}
 (-1)^{(|\lambda|+|\mu|)/2}  \det V^{p,q}_{\lambda,\mu}(\vectx;\vecta).
$$
For example, if $p=q=1$, then $F^{1,1}(\vectx ; \vecta) = a_2 - a_1$,
 and, if $p=q=2$, then $\Par_2 = \{ \emptyset, (1,1) \}$ and
\begin{align*}
F^{2,2}(\vectx ;\vecta)
 &=
\det \begin{pmatrix}
1 & x_1 & a_1 & a_1 x_1 \\
1 & x_2 & a_2 & a_2 x_2 \\
1 & x_3 & a_3 & a_3 x_3 \\
1 & x_4 & a_4 & a_4 x_4
\end{pmatrix}
-
\det \begin{pmatrix}
x_1 & x_1^2 & a_1 & a_1 x_1 \\
x_2 & x_2^2 & a_2 & a_2 x_2 \\
x_3 & x_3^2 & a_3 & a_3 x_3 \\
x_4 & x_4^2 & a_4 & a_4 x_4
\end{pmatrix}
\\
&\quad
-
\det \begin{pmatrix}
1 & x_1 & a_1 x_1 & a_1 x_1^2 \\
1 & x_2 & a_2 x_2 & a_2 x_2^2 \\
1 & x_3 & a_3 x_3 & a_3 x_3^2 \\
1 & x_4 & a_4 x_4 & a_4 x_4^2
\end{pmatrix}
+
\det \begin{pmatrix}
x_1 & x_1^2 & a_1 x_1 & a_1 x_1^2 \\
x_2 & x_2^2 & a_2 x_2 & a_2 x_2^2 \\
x_3 & x_3^2 & a_3 x_3 & a_3 x_3^2 \\
x_4 & x_4^2 & a_4 x_4 & a_4 x_4^2
\end{pmatrix}.
\end{align*}

The aim of this section is to prove the following theorem:

\begin{theorem}
\label{thm:variation}
\begin{roster}{(b)}
\item[(a)]
Let $n$ be a positive integer and let $p$ and $q$ be nonnegative integers.
For six vectors of variables
\begin{gather*}
\vectx = (x_1, \cdots, x_n),
\quad
\vecty = (y_1, \cdots, y_n),
\quad
\vectz = (z_1, \cdots, z_{p+q}),
\\
\vecta = (a_1, \cdots, a_n),
\quad
\vectb = (b_1, \cdots, b_n),
\quad
\vectc = (c_1, \cdots, c_{p+q})
\end{gather*}
we have
\begin{multline}
\det \left(
 \frac{ F^{p+1,q+1}(x_i,y_j,\vectz ; a_i,b_j,\vectc) }
      { (y_j - x_i)(1 - x_i y_j) }
 \right)_{1 \le i, j \le n}
\\
=
\frac{ (-1)^{n(n-1)/2} }{ \prod_{i,j=1}^n (y_j - x_i)(1 - x_i y_j) }
 F^{p,q}(\vectz ; \vectc)^{n-1}
 F^{n+p,n+q}(\vectx,\vecty,\vectz ; \vecta,\vectb,\vectc).
\label{variation1}
\end{multline}
\item[(b)]
Let $n$ be a positive integer and
 let $p$, $q$, $r$, $s$ be nonnegative integers.
For seven vectors of variables
\begin{gather*}
\vectx = (x_1, \cdots, x_{2n}),
\quad
\vecta = (a_1, \cdots, a_{2n}),
\quad
\vectb = (b_1, \cdots, b_{2n}),
\\
\vectz = (z_1, \cdots, z_{p+q}),
\quad
\vectc = (c_1, \cdots, c_{p+q}),
\\
\vectw = (w_1, \cdots, w_{r+s}),
\quad
\vectd = (d_1, \cdots, d_{r+s}),
\end{gather*}
we have
\begin{multline}
\Pf \left(
 \frac{ F^{p+1,q+1}(x_i,x_j,\vectz ; a_i,a_j,\vectc)
        F^{r+1,s+1}(x_i,x_j,\vectw ; b_i,b_j,\vectd) }
      { (x_j - x_i)(1 - x_i x_j) }
\right)_{1 \le i, j \le 2n}
\\
=
\frac{1}{\prod_{1 \le i < j \le 2n} (x_j - x_i)(1 - x_i x_j)}
 F^{p,q}(\vectz ; \vectc)^{n-1}
 F^{r,s}(\vectw ; \vectd)^{n-1}
\\
\times
 F^{n+p,n+q}(\vectx,\vectz ; \vecta,\vectc)
 F^{n+r,n+s}(\vectx,\vectw ; \vectb,\vectd).
\label{variation2}
\end{multline}
\end{roster}
\end{theorem}

In particular, by putting $p=q=r=s=0$ and $b_i = x_i$ for $1 \le i \le 2n$
 in (\ref{variation2}), we obtain Sundquist's identity
 \cite[Theorem~2.1]{Su}.

\begin{corollary} \label{cor:variation}
(Sundquist)
\begin{equation}
\Pf \left(
 \frac{ a_j - a_i }{ 1 - x_i x_j }
\right)_{1 \le i < j \le 2n}
 =
\frac{ (-1)^{n(n-1)/2} }
     { \prod_{1 \le i < j \le 2n} (1 - x_i x_j) }
\sum_{\lambda,\mu \in \Par_n}
 (-1)^{(|\lambda|+|\mu|)/2} \det V^{n,n}_{\lambda,\mu}(\vectx;\vecta).
\label{varition3}
\end{equation}
\end{corollary}

In order to prove Theorem~\ref{thm:variation} and
 Corollary~\ref{cor:variation}, we need a relation between
 $F^{p,q}(\vectx;\vecta)$ and $\det V^{p,q}(\vecty;\vectb)$.

\begin{prop} \label{prop:rel-fv}
We have
\begin{align}
F^{p,q}(\vectx;\vecta)
 &=
(-1)^{\binom{p}{2}+\binom{q}{2}}
\prod_{i=1}^{p+q} x_i^{p-1}
\cdot
\det V^{p,q}(\vectx+\vectx^{-1};\vecta \vectx^{q-p}),
\nonumber\\
 &=
(-1)^{\binom{p}{2}+\binom{q}{2}}
\det \U{p}{q}{\vectx}{\vectone+\vectx^2}{\vectone}{\vecta},
\label{rel-fv}
\end{align}
where
\begin{gather*}
\vectx + \vectx^{-1}
 = (x_1 + x_1^{-1}, \cdots, x_{p+q}+x_{p+q}^{-1}),
\\
\vecta \vectx^{q-p}
 = (a_1 x_1^{q-p}, \cdots, a_{p+q} x_{p+q}^{q-p}),
\\
\vectone + \vectx^2
 = (1 + x_1^{2}, \cdots, 1+x_{p+q}^{2}).
\end{gather*}
\end{prop}

Here we give a proof by using the Cauchy-Binet formula.
Let $A = \left( a_{ij} \right)_{1 \le i \le m, 1 \le j \le n}$ be
 an $m$ by $n$ matrix.
For any subsets $I = \{ i_1 < \dots < i_r \} \subset [m]$,
 and $J = \{ j_1 < \cdots < j_r \} \subset [n]$,
let $\Delta^{I}_{J}(A)$ denote the submatrix obtained by selecting the rows
 indexed by $I$ and the columns indexed by $J$.
If all rows or columns are selected, i.e., if $I = [m]$ or $J=[n]$,
 then we simply write $\Delta_J(A)$ or $\Delta^I(A)$ for
 $\Delta^{[m]}_J(A)$ or $\Delta^I_{[n]}(A)$.

\begin{lemma} \label{lem:cauchy-binet}
Let $X$ and $Y$ be any $n \times N$ matrix and $A$ be any $N \times N$ matrix.
Then we have
\begin{equation}
\det ( X A \trans Y )
 = \sum_{I,J} \det \Delta^I_J(A) \det \Delta_I(X) \det \Delta_J(Y),
\end{equation}
where the sum is taken over all pairs $(I,J)$ of $n$-element subsets
 of $[N]$.
\end{lemma}

For a partition $\lambda$ with length $\le r$, we put
$$
I(\lambda) = \{ \lambda_r, \lambda_{r-1}+1, \lambda_{r-2}+2, \cdots,
 \lambda_1+r-1 \}.
$$
Then a key of the proof of Proposition~\ref{prop:rel-fv} is
 is the following lemma.

\begin{lemma} \label{lem:minor}
Let $D_r$ be the following $r \times (2r-1)$ matrix with columns indexed by
 $0, 1, \cdots, 2r-2$:
$$
D_r =
\bordermatrix{
  & 0 &        & r-2 & r-1 & r &        & 2r-2 \cr
  &   &        &     & 1   &   &        &   \cr
  &   &        & 1   &     & 1 &        &   \cr
  &   & \rdots &     &     &   & \ddots &   \cr
  & 1 &        &     &     &   &        & 1
}.
$$
Then the minor of $D_r$ corresponding to a partition $\lambda$
 is given by
$$
\det \Delta_{I(\lambda)} \left( D_r \right)
 = \begin{cases}
(-1)^{r(r-1)/2 + |\lambda|/2} &\text{if $\lambda \in \Par_r$,} \\
0 &\text{otherwise.}
\end{cases}
$$
\end{lemma}

\begin{demo}{Proof}
First we show that $\det \Delta_{I(\lambda)} \left( D_r \right) = 0$
 unless $\lambda \in \Par_r$.
Suppose that $\det \Delta_{I(\lambda)} \left( D_r \right) \neq 0$.
Since the first row of the matrix $D_r$ has only $1$ in the $(r-1)$th column,
 we must have $r-1 \in I(\lambda)$.
If we denote by $p = p(\lambda)$ the length of the main diagonal of
 the Young diagram of $\lambda$, then we have
$$
p
 = \# \{ i : \lambda_i \ge i \}
 = \# \{ i : \lambda_i + r-i \ge r \}
 = \# \{ k \in I(\lambda) : k \ge r \}.
$$
Hence the elements $\lambda_1 + r-1 > \cdots > \lambda_p + r-p$
 are the largest $p$ elements belonging in $I(\lambda)$.
Since the Frobenius' Lemma (\cite[(1.7)]{M}) says that
$$
\{ \lambda_i + r-i : 1 \le i \le r \}
\cup
\{ r-1+j - \trans\lambda_j : 1 \le j \le r-1 \}
 =
\{ 0, 1, \cdots, 2r-2 \},
$$
we see that the elements $r - \trans\lambda_1 < \cdots < 
 r-1+p -\trans\lambda_p$ are the smallest $p$ elements not belonging
 in $I(\lambda)$.
On the other hand,
 by noting that the $k$th column of $D_r$ is identical with
 the $(2r-2-k)$th column,
 we have, if $k \neq r-1$,
 then $k \in I(\lambda)$ if and only if $2r-2-k \not\in I(\lambda)$.
Therefore we have
$$
\lambda_i + r-i + (r-1+i - \trans\lambda_i) = 2r-2
\quad(1 \le i \le p).
$$
So we have $\trans\lambda_i = \lambda_i + 1$ for $1 \le i \le p$,
 which implies $\lambda \in \Par_r$.

Next we show that, if $\lambda = (\alpha|\alpha+1) \in \Par_r$, then
 $\det\Delta_{I(\lambda)}\left(D_r\right) = (-1)^{r(r-1)/2 + |\lambda|/2}
 = (-1)^{r(r-1)/2 + |\alpha|+p}$.
Note that $\Delta_{I(\lambda)} \left( D_r \right)$ is a permutation matrix.
Let $\sigma$ be the permutation corresponding to
 $\Delta_{I(\lambda)} \left( D_r \right)$.
Then we have
\begin{gather*}
\sigma(1) > \sigma(2) > \cdots > \sigma(r-p-1),
\\
\sigma(r-p) = 1, \quad
\sigma(r-p+1) = \alpha_p+2, \quad
\cdots, \quad
\sigma(r) = \alpha_1 + 2.
\end{gather*}
The number of pairs $(i,j)$ such that $i < j $ and $\sigma(i) < \sigma(j)$
 is equal to
$$
(\alpha_1+1) + \cdots + (\alpha_p+1) = |\alpha| + p,
$$
so we have
$$
\det \Delta_{I(\lambda)}\left(D_r\right)
 = \sgn(\sigma)
 = (-1)^{r(r-1)/2-|\alpha|-p}.
$$
\end{demo}

\begin{demo}{Proof of Proposition~\ref{prop:rel-fv}}
In this proof we put $m=p+q$ for brevity.
Apply the Cauchy-Binet formula (\ref{lem:cauchy-binet}) to the following
 $(p+q) \times (2p+2q-2)$ matrices $X$ and $Y$
 (and the identity matrix $A$):
$$
X =
\begin{pmatrix}
 D_p & O \\
 O & D_q
\end{pmatrix},
\quad
Y =
\begin{pmatrix}
 1 & x_1 & \cdots & x_1^{2p-2} & a_1 & a_1 x_1 & \cdots & a_1 x_1^{2q-2} \\
 \vdots & \vdots &  & \vdots & \vdots & \vdots & & \vdots \\
 1 & x_m & \cdots & x_m^{2p-2} & a_m & a_m x_m & \cdots & a_m x_m^{2q-2}
\end{pmatrix}.
$$
Let $C \cup C' = \{ 0, 1, \cdots, 2p-2 \} \cup \{ 0', 1', \cdots, (2q-2)' \}$
 be the set indexing the columns of $X$ and $Y$.
Let $I$ be a subset of $C \cup C'$ which has cardinality $p+q$
and consider the minors
 of $X$ and $Y$ obtained by choosing the columns with indices in $I$.
Then we have
$$
\det \Delta_I\left(X\right) = 0
\quad\text{unless $\# ( I \cap C ) = p$ and $\# (I \cap C') = q$.}
$$
Suppose that $\# ( I \cap C ) = p$ and $\# (I \cap C') = q$
and that the subsets $I \cap C$ and $I \cap C'$ determine partitions
 $\lambda$ and $\mu$, i.e.,
 $I(\lambda) = I \cap C$ and $I(\mu) = I \cap C'$.
Then, by Lemma~\ref{lem:minor} and the definition of
 $V^{p,q}_{\lambda,\mu}(\vectx ; \vecta)$,
 we have
\begin{align*}
\det \Delta_I\left(X\right)
 &= \begin{cases}
 (-1)^{\binom{p}2+\binom{q}2+|\lambda|/2+|\mu|/2}
  &\text{if $\lambda \in \Par_p$ and $\mu \in \Par_q$,} \\
 0 &\text{otherwise,}
\end{cases}
\\
\det \Delta_I\left(Y\right)
 &= \det V^{p,q}_{\lambda,\mu}(\vectx;\vecta).
\end{align*}
Hence the Cauchy-Binet formula gives
\begin{align*}
\det (X \trans Y)
 &=
 \sum_{I \subset C \cup C'} \det \Delta_I\left(X\right)
 \det \Delta_I\left(Y\right)
\\
 &=
 \sum_{\lambda \in \Par_p, \, \mu \in \Par_q}
  (-1)^{\binom{p}2 + \binom{q}2 + |\lambda|/2 + |\mu|/2}
  \det V^{p,q}_{\lambda,\mu}(\vectx;\vecta).
\end{align*}

On the other hand, we have
\begin{align*}
&
\det (X \trans Y)
\\
&=
\det \begin{pmatrix}
x_1^{p-1} & x_1^p + x_1^{p-2} & \cdots & x_1^{2p-2} + 1
 & a_1 x_1^{q-1} & a_1 ( x_1^q + x_1^{q-2} ) & \cdots & a_1 ( x_1^{2q-2} + 1 )
\\
\vdots & \vdots &  & \vdots
 & \vdots & \vdots &  & \vdots
\\
x_m^{p-1} & x_m^p + x_m^{p-2} & \cdots & x_m^{2p-2} + 1
 & a_m x_m^{q-1} & a_m ( x_m^q + x_m^{q-2} ) & \cdots & a_m ( x_m^{2q-2} + 1 )
\end{pmatrix}.
\end{align*}
Now, by applying elementary transformations and by using the relations
\begin{gather*}
(x+x^{-1})^{2k}
 = \sum_{i=0}^{k-1} \binom{2k}{i} (x^{2k-2i} + x^{-2k+2i}) + \binom{2k}{k},
\\
(x+x^{-1})^{2k+1}
 = \sum_{i=0}^{k} \binom{2k+1}{i} (x^{2k-2i+1} + x^{-2k+2i-1}),
\end{gather*}
we have
\begin{align*}
\det (X \trans Y)
 &=
\prod_{i=1}^m x_i^{p-1}
\cdot
\det \left(
 \begin{cases}
  (x_i + x_i^{-1})^{j-1}
   &\text{if $1 \le i \le m$ and $1 \le j \le p$,} \\
  a_i x_i^{q-p} (x_i + x_i^{-1})^{j-p-1}
   &\text{if $1 \le i \le m$ and $p+1 \le j \le p+q$.}
 \end{cases}
\right)
\\
 &=
\prod_{i=1}^m x_i^{p-1}
 \cdot
\det V^{p,q}(\vectx+\vectx^{-1};\vecta \vectx^{q-p})
\\
 &=
\U{p}{q}{\vectx}{\vectone+\vectx^2}{\vectone}{\vecta}.
\end{align*}
\end{demo}

\begin{remark}
By the above argument in the case of $q=0$, we actually show
 one of the Littlewood's formula
$$
\sum_{\lambda \in \Par_n} s_\lambda(x_1, \cdots, x_n)
 = \prod_{1 \le i < j \le n} (1 - x_i x_j).
$$
\end{remark}

Now we are in position to derive Theorem~\ref{thm:variation}
from Theorem~\ref{thm:general}.

\begin{demo}{Proof of Theorem~\ref{thm:variation}}
Note that
\begin{equation}
\det \begin{pmatrix} x & 1 + x^2 \\ y & 1 + y^2 \end{pmatrix}
 = (x - y)(1 - xy).
\label{rel}
\end{equation}

First we prove (\ref{variation2}).
From the above relation (\ref{rel}) and (\ref{rel-fv}), we have
\begin{align*}
&
\frac{ F^{p+1,q+1}(x_i,x_j,\vectz ; a_i,a_j,\vectc)
       F^{r+1,s+1}(x_i,x_j,\vectw ; b_i,b_j,\vectd) }
     { (x_j - x_i)(1 - x_i x_j) }
\\
&\quad=
\frac{ (-1)^{1+\binom{p+1}2 + \binom{q+1}2 + \binom{r+1}2 + \binom{s+1}2} }
     { \det \begin{pmatrix}
        x_i & 1 + x_i^2 \\
        x_j & 1 + x_j^2
       \end{pmatrix}
     }
\\
&\quad\quad\times
 \det \U{p+1}{q+1}{x_i,x_j,\vectz}{1+x_i^2,1+x_j^2,\vectone+\vectz^2}
                  {1,1,\vectone}{a_i,a_j,\vectc}
\\
&\quad\quad\times
 \det \U{r+1}{s+1}{x_i,x_j,\vectw}{1+x_i^2,1+x_j^2,\vectone+\vectw^2}
                  {1,1,\vectone}{b_i,b_j,\vectd}.
\end{align*}
Hence we apply (\ref{homog2}) to obtain
\begin{align*}
&
\Pf \left(
 \frac{ F^{p+1,q+1}(x_i,x_j,\vectz ; a_i,a_j,\vectc)
        F^{r+1,s+1}(x_i,x_j,\vectw ; b_i,b_j,\vectd) }
      { (x_j - x_i)(1 - x_i x_j) }
\right)_{1 \le i, j \le 2n}
\\
&\quad=
\frac{ (-1)^{n + n \binom{p+1}2 + n \binom{q+1}2 + n \binom{r+1}2
             + n \binom{s+1}2}
     }
     { \prod_{1 \le i < j \le 2n}
         \det \begin{pmatrix}
          x_i & 1 + x_i^2 \\
          x_j & 1 + x_j^2
         \end{pmatrix}
     }
\\
&\quad\quad\times
 \det \U{p}{q}{\vectz}{\vectone+\vectz^{2}}{\vectone}{\vectc}^{n-1}
 \det \U{r}{s}{\vectw}{\vectone+\vectw^{2}}{\vectone}{\vectd}^{n-1}
\\
&\quad\quad\times
 \det \U{n+p}{n+q}{\vectx,\vectz}{\vectone+\vectx^{2},\vectone+\vectz^{2}}
                  {\vectone,\vectone}{\vecta,\vectc}
 \det \U{n+r}{n+s}{\vectx,\vectw}{\vectone+\vectx^{2},\vectone+\vectw^{2}}
                  {\vectone,\vectone}{\vectb,\vectd}
\\
&\quad=
(-1)^{n \binom{p+1}2 + n \binom{q+1}2 + n \binom{r+1}2 + n \binom{s+1}2
      +(n-1)\binom{p}{2}+ (n-1)\binom{q}{2} +(n-1)\binom{r}{2}
      +(n-1)\binom{s}{2}
      +\binom{n+p}{2} + \binom{n+q}{2}+\binom{n+r}{2} + \binom{n+s}{2}}
\\
&\quad\quad\times
\frac{ 1 }
     { \prod_{1 \le i < j \le 2n}
        (x_j - x_i)(1 - x_i x_j) }
 F^{p,q}(\vectz;\vectc)^{n-1}
 F^{r,s}(\vectw;\vectd)^{n-1}
\\
&\quad\quad\times
 F^{n+p,n+q}(\vectx,\vectz;\vecta,\vectc)
 F^{n+r,n+s}(\vectx,\vectw;\vectb,\vectd).
\end{align*}
If we use the relation
$$
n \binom{p+1}{2} - (n-1) \binom{p}{2} - \binom{p+n}{2}
 = -\binom{n}2,
$$
then we obtain the desired identity.
This proves (\ref{variation2}).

The determinant identity (\ref{variation1}) can be proven
 by the same method by using (\ref{rel}), (\ref{rel-fv}) and
 (\ref{homog1}), so we omit the detailed proof.
\end{demo}

\begin{demo}{Proof of Corollary~\ref{cor:variation}}
If we put $p=q=r=s=0$ and $b_i = x_i$ ($1 \le i \le 2n$)
 in (\ref{variation2}), then we have
$$
\Pf \left( \frac{a_j - a_i}{1 - x_i x_j} \right)_{1 \le i, j \le 2n}
 =
\frac{1}{ \prod_{1 \le i < j \le 2n} (x_j - x_i)(1 - x_i x_j) }
F^{n,n}(\vectx;\vecta) F^{n,n}(\vectx;\vectx).
$$
From the relation (\ref{rel-fv}), we see that
$$
F^{n,n}(\vectx ; \vectx)
= \prod_{i=1}^{2n} x_i^{n-1} \det V^{n,n}(\vectx+\vectx^{-1};\vectx).
$$
And, by applying appropriate elementary column transformations, we obtain
$$
F^{n,n}(\vectx;\vectx)
 =(-1)^{n(n-1)/2}
\det \left( x_i^{j-1} \right)_{1 \le i, j \le 2n}
 =(-1)^{n(n-1)/2}
\prod_{1 \le i < j \le 2n} (x_j - x_i).
$$
This completes the proof of the corollary.
\end{demo}


In Theorem~\ref{thm:variation}, we can replace
 $F^{p,q}(\vectx ; \vecta)$ by the following linear combination
 of $\det V^{p,q}_{\lambda,\mu}(\vectx ; \vecta)$:
\begin{align*}
G^{p,q}(\vectx;\vecta)
 &=
 \sum_{\lambda \in \Q_p, \, \mu \in \Q_q}
  (-1)^{(|\lambda|+|\mu|)/2} \det V^{p,q}_{\lambda,\mu}(\vectx;\vecta),
\\
H^{p,q}(\vectx;\vecta)
 &=
 \sum_{\lambda \in \R_p, \, \mu \in \R_q}
  (-1)^{(|\lambda|+p(\lambda)+|\mu|+p(\mu))/2}
  \det V^{p,q}_{\lambda,\mu}(\vectx;\vecta),
\end{align*}
where $\Q_n$ (resp. $\R_n$) is the set of partitions $\lambda$
 with length $\le n$ which is of the form $\lambda = (\alpha+1|\alpha)$
 (resp. $\lambda = (\alpha|\alpha)$) in the Frobenius notation.

The following Lemma and Proposition can be proven by the same
 idea as Lemma~\ref{lem:minor} and Proposition~\ref{prop:rel-fv},
so we leave the proof to the reader.

\begin{lemma}
Let $C_r$ (resp. $B_r$) be the $r \times (2r+1)$
 (resp. $r \times 2r$) matrix given by
\begin{gather*}
B_r =
\bordermatrix{
  & 0  &        & r-2 & r-1 & r & r+1 &        & 2r-1 \cr
  &    &        &     & -1  & 1 &     &        &   \cr
  &    &        & -1  &     &   & 1   &        &   \cr
  &    & \rdots &     &     &   &     & \ddots &   \cr
  & -1 &        &     &     &   &     &        & 1
},
\\
C_r =
\bordermatrix{
  & 0  &        & r-2 & r-1 & r & r+1 & r+2 &        & 2r \cr
  &    &        &     & -1  & 0 & 1   &     &        &   \cr
  &    &        & -1  &     &   &     & 1   &        &   \cr
  &    & \rdots &     &     &   &     &     & \ddots &   \cr
  & -1 &        &     &     &   &     &     &        & 1
}.
\end{gather*}
Then we have
\begin{roster}{(2)}
\item[(1)]
For a partition $\lambda$ of length $\le r$, we have
$$
\det \Delta_{I(\lambda)} \left( B_r \right)
 = \begin{cases}
 (-1)^{\binom{r+1}{2} + (|\lambda|+p(\lambda))/2}
 &\text{if $\lambda \in \R_r$,} \\
 0 &\text{otherwise.}
\end{cases}
$$
\item[(2)]
For a partition $\lambda$ of length $\le r$, we have
$$
\det \Delta_{I(\lambda)} \left( C_r \right)
 = \begin{cases}
 (-1)^{\binom{r+1}{2} + |\lambda|/2}
 &\text{if $\lambda \in \Q_r$,} \\
 0 &\text{otherwise.}
\end{cases}
$$
\end{roster}
\end{lemma}

\begin{prop}
\begin{align*}
G^{p,q}(\vectx;\vecta)
 &=
 (-1)^{\binom{p}2+\binom{q}2}
 \prod_{i=1}^{p+q} x_i^{p-1} (1 - x_i^2) \cdot
 \det V^{p,q}(\vectx+\vectx^{-1};\vecta \vectx^{q-p}),
\\
 &=
 (-1)^{\binom{p}2+\binom{q}2}
 \prod_{i=1}^{p+q} (1 - x_i^2) \cdot
 \det \U{p}{q}{\vectx}{\vectone+\vectx^{2}}{\vectone}{\vecta},
\\
H^{p,q}(\vectx;\vecta)
 &=
 (-1)^{\binom{p}2+\binom{q}2}
 \prod_{i=1}^{p+q} x_i^{p-1} (1 - x_i) \cdot
 \det V^{p,q}(\vectx+\vectx^{-1};\vecta \vectx^{q-p}),
\\
 &=
 (-1)^{\binom{p}2+\binom{q}2}
 \prod_{i=1}^{p+q} (1 - x_i) \cdot
 \det \U{p}{q}{\vectx}{\vectone+\vectx^{2}}{\vectone}{\vecta}.
\end{align*}
In particular, we have
\begin{align*}
G^{p,q}(\vectx;\vecta)
 &= \prod_{i=1}^{p+q} (1 - x_i^2) \cdot F^{p,q}(\vectx;\vecta),
\\
H^{p,q}(\vectx;\vecta)
 &= \prod_{i=1}^{p+q} (1 - x_i) \cdot F^{p,q}(\vectx;\vecta).
\end{align*}
\end{prop}

From these relations, we have determinant and Pfaffian identities
 involving $G^{p,q}(\vectx;\vecta)$ and $H^{p,q}(\vectx;\vecta)$
 similar to (\ref{variation1}) and (\ref{variation2}).
More generally, we can consider, for example,
$$
\sum_{\lambda \in \Par_p, \, \mu \in \Q_q}
 (-1)^{(|\lambda|+|\mu|)/2}
 \det V^{p,q}_{\lambda,\mu}(\vectx;\vecta),
$$
which can be expressed in terms of $\det V^{p,q}$ or $\det U^{p,q}$.


\section{
Another generalization of Cauchy's determinant identity
}

In this section, we give another type of generalized Cauchy's determinant
 identities involving $\det V^{p,q}$ and $\det W^p$.

\begin{theorem} \label{thm:another}
\begin{align}
&
\det \left(
 \frac{ 1 }
      { \det V^{p+1,q+1}(x_i,y_j,\vectz ; a_i, b_j, \vectc) }
\right)_{1 \le i, j \le n}
\notag
\\
&\quad=
\frac{ (-1)^{n(n-1)/2}
       \prod_{1 \le i < j \le n}
         \det V^{p+1,q+1}(x_i,x_j,\vectz ; a_i,a_j,\vectc)
         \det V^{p+1,q+1}(y_i,y_j,\vectz ; b_i,b_j,\vectc) }
     { \prod_{i,j=1}^n \det V^{p+1,q+1}(x_i,y_j,\vectz ; a_i,b_j,\vectc) },
\label{another1}
\\
&
\det \left(
 \frac{ 1 }
      { \det W^{p+2}(x_i,y_j,\vectz ; a_i, b_j, \vectc) }
\right)_{1 \le i, j \le n}
\notag
\\
&\quad=
\frac{ (-1)^{n(n-1)/2}
       \prod_{1 \le i < j \le n}
         \det W^{p+2}(x_i,x_j,\vectz ; a_i,a_j,\vectc)
         \det W^{p+2}(y_i,y_j,\vectz ; b_i,b_j,\vectc) }
     { \prod_{i,j=1}^n \det W^{p+2}(x_i,y_j,\vectz ; a_i,b_j,\vectc) }.
\label{another2}
\end{align}
\end{theorem}

If $p=q=0$, then the identity (\ref{another1}) becomes
$$
\det \left( \frac{1}{b_j - a_i} \right)_{1 \le i, j \le n}
 =
\frac{ (-1)^{n(n-1)/2} \prod_{1 \le i < j \le n} (a_j - a_i)(b_j - b_i) }
     { \prod_{i,j=1}^n (b_j - a_i) },
$$
which is equivalent to Cauchy's determinant identity (\ref{cauchy}).

In order to prove Theorem~\ref{thm:another}, we put
$$
f(x, y ; a, b)
 =
 \det V^{p+1,q+1}(x,y,\vectz;a,b,\vectc),
\quad\text{or}\quad
 \det W^{p+2}(x,y,\vectz;a,b,\vectc).
$$
The proof is based on the quadratic relations among $f(x,y;a,b)$'s,
 which follows from the Pl\"ucker relation for determinants.

\begin{lemma} \label{lem:plucker}
Let $\left( a_1, a_2, a_3, a_4, x_1, \cdots, x_m \right)$ be
 a $(m+2) \times (m+4)$ matrix.
If we put
$$
D(i,j)
 =
 \det \left( a_i, a_j, x_1, \cdots, x_m \right),
$$
then we have
$$
D(1,2) D(3,4) - D(1,3) D(2,4) + D(1,4) D(2,3) = 0.
$$
\end{lemma}

\begin{prop}
The polynomials $f(x,y ; a,b)$'s satisfy
\begin{multline}
f(x_1, x_2 ; a_1, a_2) f(y_1, y_2 ; b_1, b_2)
-
f(x_1, y_1 ; a_1, b_1) f(x_2, y_2 ; a_2, b_2)
+
f(x_1, y_2 ; a_1, b_2) f(x_2, y_1 ; a_2, b_1)
\\
=
0.
\label{plucker-vw}
\end{multline}
\end{prop}

\begin{demo}{Proof}
Apply the Pl\"ucker relation to the transposes of the matrices
$$
\begin{pmatrix}
1 & x_1 & \cdots & a_1 x_1^q \\
1 & x_2 & \cdots & a_2 x_2^q \\
1 & y_1 & \cdots & b_1 y_1^q \\
1 & y_2 & \cdots & b_2 y_2^q \\
1 & z_1 & \cdots & c_1 z_1^q \\
\vdots & \vdots & & \vdots   \\
1 & z_m & \cdots & c_m z_m^q \\
\end{pmatrix}
\quad\text{or}\quad
\begin{pmatrix}
1 + a_1 x_1^{p+1} & x_1 + a_1 x_1^p & \cdots & x_1^{p+1} + a_1 \\
1 + a_2 x_2^{p+1} & x_2 + a_2 x_2^p & \cdots & x_2^{p+1} + a_2 \\
1 + b_1 y_1^{p+1} & y_1 + b_1 y_1^p & \cdots & y_1^{p+1} + b_1 \\
1 + b_2 y_2^{p+1} & y_2 + b_2 y_2^p & \cdots & y_2^{p+1} + b_2 \\
1 + c_1 z_1^{p+1} & z_1 + c_1 z_1^p & \cdots & z_1^{p+1} + c_1 \\
\vdots            & \vdots          &        & \vdots          \\
1 + c_p z_p^{p+1} & z_p + c_p z_p^p & \cdots & z_p^{p+1} + c_p
\end{pmatrix},
$$
where $m=p+q$.
\end{demo}

\begin{demo}{Proof of Theorem~\ref{thm:another}}
We proceed by induction on $n$.
In this proof, we write $f(x_i,x_j)$ (resp. $f(x_i,y_j)$, $f(y_i,y_j)$)
 instead of $f(x_i,x_j;a_i,a_j)$ (resp. $f(x_i,y_j;a_i,b_j)$,
 $f(y_i,y_j;b_i,b_j)$).

If $n=1$, then there is nothing to prove, and,
 if $n=2$, then the desired identities are equivalent to
 the quadratic relations (\ref{plucker-vw}).

Suppose $n \ge 3$.
Applying the Desnanot--Jacobi formula (\ref{det-dodgson}) to
 the matrix
$$
A = \left( \frac{1}{ f(x_i,y_j) } \right)_{1 \le i, j \le n}.
$$
Then, from the induction hypothesis, we have
$$
\det A^k_l
 =
\frac{ (-1)^{(n-1)(n-2)/2} 
       \prod_{i=3}^n f(x_{k'},x_i) f(y_{l'},y_i)
       \prod_{3 \le i < j \le n} f(x_i,x_j) f(y_i,y_j)
     }
     { f(x_{k'},y_{l'}) \prod_{i=3}^n f(x_{k'},y_i) f(x_i,y_{l'})
       \prod_{i,j=3}^n f(x_i,y_j)
     }
$$
where $1 \le k, l \le 2$ and the indices $k'$ and $l'$ are determined by
 the condition $\{ k, k' \} = \{ l, l' \} = \{ 1, 2 \}$,
 and
$$
\det A^{1,2}_{1,2}
 =
\frac{ (-1)^{(n-2)(n-3)/2}
       \prod_{3 \le i < j \le n} f(x_i,x_j)
       \prod_{3 \le i < j \le n} f_(y_i,y_j) }
     { \prod_{i,j=3}^n f(x_i,y_j) }.
$$
By cancelling the common factors, we see that, in order to prove
 the identity, it is enough to show
$$
- f(x_1,x_2) f(y_1,y_2) = f(x_1,y_2) f(x_2,y_1) - f(x_1,y_1) f(x_2,y_2).
$$
This is equivalent to (\ref{plucker-vw}).
\end{demo}

\section{A hyperpfaffian expression}

H.~Tagawa finds that $\det V^{n,n}(\vectx;\vecta)$ is expressed
 by a hyperpfaffian.
The aim of this section is to prove this expression.

First we recall the definition of hyperpfaffians (see \cite{LT}).
Let $n$ and $r$ be positive integers.
Define a subset ${\mathcal{E}}_{rn,n}$ of the symmetric groups
 $\Sym_{rn}$ by
$$
\mathcal{E}_{rn,n}
 =
\left\{
 \sigma \in \Sym_{rn}
 :
 \sigma(n(i-1)+1) < \sigma(n(i-1)+2) < \cdots < \sigma(ni)
 \text{ for $1 \le i \le r$}
\right\}.
$$
For example, if $n=r=2$, then $\mathcal{E}_{4,2}$ is composed
 of the following 6 elements:
$$
\mathcal{E}_{4,2}
 =
\left\{
 (1,2,3,4), (1,3,2,4), (1,4,2,3),
 (3,4,1,2), (2,4,1,3), (2,3,1,4)
\right\}.
$$
Let $a = \left( a_{i_1 \dots i_n} \right)_{1 \le i_1 < \dots < i_n \le nr}$
 be an alternating tensor,
i.e. $a_{i_{\sigma(1)} \dots i_{\sigma(n)}} = \sgn(\sigma) a_{i_1 \dots i_n}$
 for any permutations $\sigma \in \Sym_{nr}$.
The hyperpfaffian of $a$ is, by definition,
$$
\Pf^{[n]}(a)
 =
\frac{1}{r!}
\sum_{\sigma \in \mathcal{E}_{nr,n}}
 \sgn(\sigma)
 \prod_{i=1}^{r} a_{\sigma(n(i-1)+1),\dots, \sigma(ni)}.
$$

An alternating $2$ tensor $a$ is a skew-symmetric matrix
 and the hyperpfaffian $\Pf^{[2]}(a)$ is the usual Pfaffian of
 the skew-symmetric matrix.
J.-G.~Luque and J.-Y.~Thibon \cite{LT} computed the following composition
 of hyperpfaffians by using the Grassmann algebra.

\begin{prop} (\cite{LT})
\label{prop:compo}
Let $n$ and $r$ be positive integers and assume $n=2m$ is even.
Given a skew-symmetric matrix $A = (a_{ij})_{1 \le i,j \le nr}$,
 we define an alternating $n$-tensor $A^{[n]}$ by putting
$$
\left( A^{[n]} \right)_{i_1, \cdots, i_n}
 = \Pf \left( a_{i_k,i_l} \right)_{1 \le k,l \le n}
\quad\text{for $1 \le i_1 < \cdots < i_n \le nr$}.
$$
Then we have
\begin{equation}
\Pf^{[n]}(A^{[n]})
= \frac{(mr)!}{(m!)^r r!} \Pf (A).
\label{compo}
\end{equation}
\end{prop}

The main theorem in this section is the following.

\begin{theorem} \label{thm:hyper}
If $n$ is even,
then
\begin{gather}
\det V^{n,n}(\vectx ; \vecta)
 =
\Pf^{[n]} \left[
 \left( 1 + \prod_{s=1}^n a_{i_s} \right)
 \prod_{1 \le s < t \le n} (x_{i_t} - x_{i_s})
\right]_{1 \le i_1 < \dots < i_n \le 2n},
\label{hyper-v}
\\
\det \U{n}{n}{\vectx}{\vecty}{\vecta}{\vectb}
 =
\Pf^{[n]} \left[
 \left( \prod_{s=1}^n a_{i_s} + \prod_{s=1}^n b_{i_s} \right)
 \prod_{1 \le s < t \le n}
  \det \begin{pmatrix}
   y_{i_s} & x_{i_s} \\
   y_{i_t} & x_{i_t}
  \end{pmatrix}
\right]_{1 \le i_1 < \dots < i_n \le 2n}.
\label{hyper-u}
\end{gather}
\end{theorem}

To prove this theorem, we need to compute the following special
 Pfaffian and hyperpfaffian.

\begin{lemma} \label{lem:special-pf}
Let $n$ and $r$ be positive integers and assume $n = 2m$ is even.
Then we have
\begin{gather}
\Pf \left( \frac{ (x_j^m - x_i^m)^2 }{ x_j - x_i } \right)_{1 \le i, j \le nr}
 =
\begin{cases}
  \prod_{1 \le i < j \le n} (x_j - x_i) &\text{if $r = 1$,} \\
  0 &\text{if $r \ge 2$,}
\end{cases}
\label{special-pf}
\\
\Pf^{[n]} \left[
 \prod_{1 \le s < t \le n} (x_{i_t}-x_{i_s})
\right]_{1 \le i_1 < \dots < i_n \le nr}
 =
\begin{cases}
  \prod_{1 \le i < j \le n} (x_j - x_i) &\text{if $r = 1$,} \\
  0 &\text{if $r \ge 2$.}
\end{cases}
\label{special-hyppf}
\end{gather}
\end{lemma}

\begin{demo}{Proof}
There are several ways to prove the identity (\ref{special-pf}).
Here we appeal to Theorem~\ref{thm:main} (\ref{main2}).
If we take $p=q=r=s=0$ and put $a_i = b_i = x_i^m$ ($1 \le i \le nr$) in
 (\ref{main2}), we have
$$
\Pf \left( \frac{ (x_j^m - x_i^m)^2 }{ x_j - x_i } \right)_{1 \le i, j \le nr}
 =
\frac{1}{ \prod_{1 \le i < j \le nr} (x_j - x_i) }
\det V^{mr,mr}(\vectx ; \vectx^m)^2.
$$
If $r = 1$, then $\det V^{m,m}(\vectx ; \vectx^m)$ is the usual Vandermonde
 determinant and $\det V^{m,m}(\vectx ; \vectx^m)
 = \prod_{1 \le i < j \le n}(x_j - x_i)$.
If $r \ge 2$, then the $(m+1)$st column of $V^{mr,mr}(\vectx;\vectx^m)$
 is the same as the $(rm+1)$st column, so we have
 $\det V^{mr,mr}(\vectx;\vectx^m) = 0$.
Hence we obtain (\ref{special-pf}).

Next we prove (\ref{special-hyppf}).
Apply Proposition~\ref{prop:compo} to the matrix
 $A = \left( (x_j^m - x_i^m)^2/(x_j - x_i) \right)_{1 \le i, j \le nr}$.
Then it follows from (\ref{special-pf}) that
\begin{align*}
\Pf^{[n]} \left[
 \prod_{1 \le s < t \le n} (x_{i_t} - x_{i_s})
\right]_{1 \le i_1 < \cdots < i_n \le nr}
 &=
\Pf^{[n]} \left[
 \Pf \left(
   \frac{ ( x_{i_l}^m - x_{i_k}^m )^2 }{ x_{i_l} - x_{i_k} }
 \right)_{1 \le k, l \le n}
\right]_{1 \le i_1 < \cdots < i_n \le nr}
\\
 &=
\frac{(mr)!}{(m!)^r r!}
\Pf \left(
 \frac{ ( x_j^m - x_i^m )^2 }{ x_j - x_i }
\right)_{1 \le i, j \le nr}.
\end{align*}
Again using (\ref{special-pf}), we obtain the desired identity.
\end{demo}

Now we are in position to prove Theorem~\ref{thm:hyper}.

\begin{demo}{Proof of Theorem~\ref{thm:hyper}}
The identity (\ref{hyper-u}) immediately follows from (\ref{hyper-v})
 by noting the relation (\ref{rel-uv2}),
 so we prove (\ref{hyper-v}).

Let $\binom{[2n]}{n}$ denote the set of all $n$-element subsets of $[2n]
 = \{ 1, 2, \cdots, 2n \}$.
For a subset $I \in \binom{[2n]}{n}$, we put
$$
a_I = \prod_{i \in I} a_i,
\quad
\Delta(\vectx_I) = \prod_{\substack{i,j \in I \\ i < j}} (x_j - x_i).
$$
By the Laplace expansion formula and Vandermonde determinant formula,
 we have
$$
V^{n,n}(\vectx;\vecta)
 =
\sum_{I \in \binom{[2n]}{n}} (-1)^{|I|+\binom{n+1}{2}}
 a_{I} \Delta(\vectx_I) \Delta(\vectx_{I^c}),
$$
where we write $|I| = \sum_{i \in I} i$ and denote by $I^c$
 the complementary subset $I$ in $[2n]$.

On the other hand, by the definition of hyperpfaffians,
we see that
\begin{align*}
\Pf^{[n]} \left[
 ( 1 + a_I ) \Delta(\vectx_I)
\right]_I
 &=
\frac{1}{2!}
\sum_{I \in \binom{[2n]}{n}} (-1)^{|I|+\binom{n+1}{2}}
 ( 1 + a_I ) \Delta(\vectx_I) ( 1 + a_{I^c} ) \Delta(\vectx_{I^c})
\\
 &=
\left( 1 + \prod_{i=1}^{2n} a_i \right)
\Pf^{[n]} \left[
 \prod_{1 \le s < t \le n} \left( x_{i_t}-x_{i_s} \right)
\right]_{1 \le i_1 < \dots < i_n \le 2n}
\\
 &\quad+
\sum_{I \in \binom{[2n]}{n}}
 (-1)^{|I|+\binom{n+1}{2}}
 a_I \Delta(\vectx_I) \Delta(\vectx_{I^c}).
\end{align*}
By (\ref{special-hyppf}), the first term vanishes,
and we obtain the desired formula.
\end{demo}
At the end of this section,
we should remark that H.~Tagawa has a similar hyperpfaffian expression
for the case that $p=q$ is odd.
It has slightly different from the case $p=q$ is even,
but we don't have any general formula when $p\neq q$.

\section{
Application to Littlewood--Richardson coefficients
}

In this section, we use the Pfaffian identity (\ref{main3}) in Theorem 1.1
 and the minor-summation formula \cite{IW1}
 to derive a relation between Littlewood--Richardson coefficients.

For three partitions $\lambda$, $\mu$ and $\nu$,
 we denote by $\LR^\lambda_{\mu,\nu}$ the Littlewood-Richardson coefficient.
These numbers $\LR^\lambda_{\mu,\nu}$ appear in the following expansions
 (see \cite{M}) :
\begin{gather*}
s_\mu(X) s_\nu(X)
 = \sum_\lambda \LR^\lambda_{\mu,\nu} s_\lambda(X),
\\
s_{\lambda/\mu}(X)
 = \sum_\nu \LR^\lambda_{\mu,\nu} s_\nu(X),
\\
s_\lambda(X, Y)
 = \sum_{\mu,\nu} \LR^\lambda_{\mu,\nu} s_\mu(X) s_\nu(Y).
\end{gather*}

We are concerned with the Littlewood--Richardson coefficients
 involving rectangular partitions.
Let $\square(a,b)$ denote the partition whose Young diagram is
 the rectangle $a\times b$,
i.e.
$$
\square(a,b) = (b^a) = (\underbrace{b,\dots,b}_{a}).
$$
For a partition $\lambda \subset \square(a,b)$, we define a partition
 $\lambda^\dagger = \lambda^\dagger(a,b)$ by
$$
\lambda^\dagger_i = b - \lambda_{a+1-i}
\quad(1 \le i \le a).
$$
This partition $\lambda^\dagger$ is the complement of $\lambda$
 in the rectangle $\square(a,b)$.

Okada \cite{O1} used the special case of the identities
 (\ref{main1}) and (\ref{main2}) (i.e., the case of $p=q=0$
 and $p=q=r=s=0$) to prove the following proposition.

\begin{prop} \label{prop:LR}
Let $n$ be a positive integer and let $e$ and $f$ be nonnegative integers.
\begin{roster}{(2)}
\item[(1)]
For partitions $\mu$, $\nu$, we have
\begin{equation}
\LR^{\square(n,e)}_{\mu,\nu}
 = \begin{cases}
 1 &\text{if $\nu = \mu^\dagger(n,e)$,} \\
 0 &\text{otherwise.}
\end{cases}
\label{LR1}
\end{equation}
\item[(2)]
For a partition $\lambda$ of length $\le 2n$, we have
\begin{equation}
\LR^\lambda_{\square(n,e), \square(n,f)}
 = \begin{cases}
 1 &\text{if $\lambda_{n+1} \le \min(e,f)$
 and $\lambda_i + \lambda_{2n+1-i} = e+f$ ($1 \le i \le n$),}
 \\
 0 &\text{otherwise.}
\end{cases}
\label{LR2}
\end{equation}
\end{roster}
\end{prop}

The main result of this section is the following theorem,
 which generalizes (\ref{LR2}).

\begin{theorem} \label{thm:LR}
Let $n$ be a positive integer and let $e$ and $f$ be nonnegative integers.
Let $\lambda$ and $\mu$ be partitions such that the length $l(\lambda) \le 2n$
 and $\mu \subset \square(n,e)$.
Then we have
\begin{roster}{(2)}
\item[(1)]
$\LR^\lambda_{\mu,\square(n,f)} = 0$ unless
\begin{equation}
\lambda_n \ge f
\quad\text{and}
\quad
\lambda_{n+1} \le \min (e,f).
\label{cond}
\end{equation}
\item[(2)]
If $\lambda$ satisfies the above condition (\ref{cond}) and
 we define two partitions $\alpha$ and $\beta$ by
\begin{equation}
\alpha_i = \lambda_i - f,
\quad
\beta_i = e - \lambda_{2n+1-i},
\quad
(1 \le i \le n),
\label{def-ab}
\end{equation}
then we have
$$
\LR^\lambda_{\mu,\square(n,f)}
 = \LR^\beta_{\alpha, \mu^\dagger(n,e)}.
$$
In particular, $\LR^\lambda_{\mu,\square(n,f)} = 0$ unless $\alpha \subset
 \beta$.
\end{roster}
\end{theorem}

In particular, if $\mu = \square(n,e)$ is a rectangle,
 then this theorem reduces to (\ref{LR2}),
 because $\LR^\alpha_{\beta,\emptyset} = \delta_{\alpha,\beta}$.
If $\mu$ is a near-rectangle, then we have the following corollary
 by using Pieri's rule \cite[(5.16), (5.17)]{M}.

\begin{corollary}
Suppose that a partition $\lambda \subset \square(2n,e+f)$ satisfies
 the condition (\ref{cond}) in Theorem~\ref{thm:LR}.
Define two partitions $\alpha$ and $\beta$ by (\ref{def-ab}).
Then we have
\begin{align*}
\LR^\lambda_{(e^{n-1},e-k), (f^n)}
&= \begin{cases}
1 &\text{if $\beta/\alpha$ is a horizontal strip of length $k$,} \\
0 &\text{otherwise,}
\end{cases}
\\
\LR^\lambda_{(e^{n-k},(e-1)^k), (f^n)}
&= \begin{cases}
1 &\text{if $\beta/\alpha$ is a vertical strip of length $k$,} \\
0 &\text{otherwise.}
\end{cases}
\end{align*}
\end{corollary}

In order to prove Theorem~\ref{thm:LR}, we substitute
\begin{equation}
a_i = x_i^{e+p+n},
\quad
c_i = z_i^{e+p+n},
\quad
d_i = w_i^{f+r+n}
\label{subs4}
\end{equation}
in the Pfaffian identity (\ref{main2}).
By the bi-determinant definition of Schur functions, we have
$$
\det V^{p,q}(\vectx ; \vectx^k)
 = \begin{cases}
 s_{\square(q,k-p)}(\vectx) \Delta(\vectx)
 &\text{if $k \ge p$,} \\
 0 &\text{if $k < p$,}
\end{cases}
$$
where $\Delta(\vectx) = \prod_{1 \le i < j \le n} (x_j - x_i)$.
Hence, under the substitution (\ref{subs4}), the identity (\ref{main2})
 give us the following Pfaffian identity.

\begin{prop}
We have
\begin{multline}
\frac{1}{\Delta(\vectx)}
 \Pf \left(
  (x_j - x_i)
  s_{\square(q+1,e+n-1)}(x_i,x_j,\vectz)
  s_{\square(s+1,f+n-1)}(x_i,x_j,\vectw)
 \right)_{1 \le i, j \le 2n}
\\
=
s_{\square(q,e+n)}(\vectz)^{n-1}
s_{\square(s,f+n)}(\vectw)^{n-1}
s_{\square(n+q,e)}(\vectx, \vectz)
s_{\square(n+s,f)}(\vectx, \vectw).
\label{pf-schur}
\end{multline}
\end{prop}

\begin{remark}
If we substitute
$$
a_i = x_i^{e+p+n},
\quad
b_i = y_i^{e+p+n}
\quad
(1 \le i \le n)
$$
in the determinant identity (\ref{main1}), then we have
\begin{multline}
\frac{1}{\Delta(\vectx) \Delta(\vecty)}
 \det \left(
  s_{\square(q+1,e+n-1)}(x_i,y_j,\vectz)
 \right)_{1 \le i, j \le n}
\\
=
(-1)^{n(n-1)/2}
s_{\square(q,e+n)}(\vectz)^{n-1}
s_{\square(q+n,e)}(\vectx,\vecty,\vectz).
\label{det-schur}
\end{multline}
The special case ($q=e+n-1$) of this identity is given in
 \cite[Proposition~8.4.3]{L}, and the proof there works
 in the general case.
\end{remark}

If we take $q=s=0$ in (\ref{pf-schur}), we have
\begin{multline}
\frac{1}{\Delta(\vectx)}
\Pf \left(
(x_j - x_i)
h_{e+n-1}(x_i,x_j,\vectz)
h_{f+n-1}(x_i,x_j,\vectw)
\right)_{1 \le i, j \le 2n}
\\
=
s_{\square(n,e)}(\vectx, \vectz)
s_{\square(n,f)}(\vectx, \vectw).
\label{pf-schur2}
\end{multline}
We use the minor-summation formula \cite{IW1} to expand the left hand side
 in the Schur function bases $\{ s_\lambda( \vectx ) \}$.

\begin{lemma} \label{lem:coeff}
Let $b_{k,l}$ be the coefficient of $x^k y^l$ in
$$
(y - x) h_{e+n-1}(x,y,\vectz) h_{f+n-1}(x,y,\vectw).
$$
Then we have $b_{kl} = - b_{lk}$, and
 $b_{kl}$, $k < l$, is given by
$$
b_{kl}
 = \sum_{i,j} h_i(\vectz) h_j(\vectw),
$$
where the sum is taken over all pairs of integers $(i,j)$ satisfying
$$
i+j = (e+n-1)+(f+n-1)+1-k-l,
\quad
0 \le i \le (e+n-1)-k,
\quad
0 \le j \le (f+n-1)-k.
$$
\end{lemma}

Note that $b_{kl} = 0$ unless $0 \le k, l \le e+f+2n-1$.

\begin{demo}{Proof}
By using the relation
$$
h_r(x,y,\vectz)
 = \sum_{a,b \ge 0} x^a y^b h_{r-a-b}(\vectz),
$$
we see that
$$
b_{kl}
 =
\left(
 \sum_{\substack{
 0 \le a, b \le e+n-1, \ 0 \le c, d \le f+n-1 \\ a + c = k, \ b + d = l-1
 }}
-
\sum_{\substack{
 0 \le a, b \le e+n-1, \ 0 \le c, d \le f+n-1 \\ a + c = k-1, \ b + d = l
 }}
\right)
h_{(e+n-1)-a-b}(\vectz) h_{(f+n-1)-c-d}(\vectw).
$$
Let $b_{kl}(i,j)$ be the coefficient of $h_i(\vectz) h_j(\vectw)$ in $b_{kl}$.
Then, by considering the homogeneous degree, we see that $b_{kl}(i,j) = 0$
 unless $i + j = (e+n-1) + (f+n-1) + 1 - k - l$, $0 \le i \le e+n-1$ and
 $0 \le j \le f+n-1$.

Now we assume
$$
i + j = (e+n-1) + (f+n-1) + 1 - k - l,
\quad
0 \le i \le e+n-1,
\quad
0 \le j \le f+n-1.
$$
If we put
\begin{align*}
S_1
&=
\left\{
 \begin{pmatrix} a \\ b \\ c \\ d \end{pmatrix} \in \Nat^4 :
 \begin{matrix}
  a + c = k, \\ b + d = l-1, \\ a + b = (e+n-1)-i, \\ c + d = (f+n-1)-j
 \end{matrix}
\right\},
\\
S_2
&=
\left\{
 \begin{pmatrix} a' \\ b' \\ c' \\ d' \end{pmatrix} \in \Nat^4 :
 \begin{matrix}
  a' + c' = k-1, \\ b' + d' = l, \\ a' + b' = (e+n-1)-i, \\ c' + d' = (f+n-1)-j
 \end{matrix}
\right\},
\end{align*}
where $\Nat$ denotes the set of nonnegative integers, then we have
$$
b_{kl}(i,j) = \# S_1 - \# S_2.
$$
The solutions to the equations in $S_1$ and $S_2$ are given by
$$
\begin{pmatrix} a \\ b \\ c \\ d \end{pmatrix}
 =
\begin{pmatrix}
 t + (e+n-1)-i-l+1 \\ -t + l-1 \\ -t + k+l-1-(e+n-1)+i \\ t
\end{pmatrix}
,
\quad\quad
\begin{pmatrix} a' \\ b' \\ c' \\ d' \end{pmatrix}
 =
\begin{pmatrix}
 t + (e+n-1)-i-l \\ -t + l \\ -t + k+l-1-(e+n-1)+i \\ t
\end{pmatrix}.
$$
Hence we see that
\begin{align*}
\# S_1
 &= \# \{ t \in \Int : t \ge a_0, \ t \le b_0, \ t \le c_0, \ t \ge d_0 \},
\\
\# S_2
 &= \# \{ t \in \Int : t \ge a'_0, \ t \le b'_0, \ t \le c'_0, \ t \ge d'_0 \},
\end{align*}
where
$$
\begin{pmatrix} a_0 \\ b_0 \\ c_0 \\ d_0 \end{pmatrix}
=
\begin{pmatrix}
l-1+i-(e+n-1) \\ l-1 \\ k+l-1-(e+n-1)+i \\ 0 \end{pmatrix},
\quad\quad
\begin{pmatrix} a'_0 \\ b'_0 \\ c'_0 \\ d'_0 \end{pmatrix}
=
\begin{pmatrix}
l+i-(e+n-1) \\ l \\ k+l-1-(e+n-1)+i \\ 0 \end{pmatrix}.
$$
We compute $\# S_1$ and $\# S_2$ in the following four cases :
\begin{roster}{(c)}
\item[(a)]
$i \le (e+n-1)-k$ and $j \le (f+n-1)-k$.
\item[(b)]
$i \le (e+n-1)-k$ and $j > (f+n-1)-k$.
\item[(c)]
$i > (e+n-1)-k$ and $j \le (f+n-1)-k$.
\item[(d)]
$i > (e+n-1)-k$ and $j > (f+n-1)-k$.
\end{roster}

Here we note that
$$
j \le (f+n-1)-k \quad\text{if and only if}\quad l+i-(e+n-1)-1 \ge 0,
$$
and that
\begin{gather*}
a_0 - d_0 = l+i-(e+n-1)-1,
\quad
b_0 - c_0 = (e+n-1)-i-k,
\\
a'_0 - d'_0 = l+i-(e+n-1) = a_0 - d_0 +1,
\quad
b'_0 - c'_0 = (e+n-1)-i-k+1 = b_0 - c_0 + 1.
\end{gather*}
Hence we see that,
if $i \le (e+n-1)-k$, then $b_0 \ge c_0$ and $b'_0 > c'_0$, and that,
if $i \le l-k$, then $a_0 \ge d_0$ and $a'_0 > d'_0$.

In Case (a), we have $a_0 \ge d_0$, $b_0 \ge c_0$, $a'_0 > d'_0$
 and $b'_0 > c'_0$, so
\begin{align*}
\# S_1
 &= \# \{ t \in \Int : t \ge a_0, \ t \le c_0 \}
 = c_0 - a_0 + 1
 = k+1,
\\
\# S_2
 &= \# \{ t \in \Int : t \ge a'_0 , \ t \le c'_0 \}
 = c'_0 - a'_0 +1
 = k
\end{align*}
Hence we have $b_{kl}(p,q) = 1$.
(This argument holds if $k = 0$.)
In Case (b), we have $b_0 \ge c_0$, $b'_0 > c'_0$, $a_0 < d_0$
 and $a'_0 \le d'_0$, so
\begin{align*}
\# S_1
 &= \# \{ t \in \Int : t \ge d_0, \ t \le c_0 \}
 = c_0 - d_0 + 1
 = k+l+i-(e+n-1),
\\
\# S_2
 &= \# \{ t \in \Int : t \ge d'_0 , \ t \le c'_0 \}
 = c'_0 - d'_0 +1
 = k+l+i-(e+n-1).
\end{align*}
Hence we have $b_{kl}(i,j) = 0$.
Similarly, in Case (c), we have $b_{kl}(i,j) = 0$.
In Case (d), we have $i+j > (e+n-1)+(f+n-1)-k-l+1$,
 which contradicts to the assumption $i+j = (e+n-1)+(f+n-1)-k-l+1$.

This completes the proof.
\end{demo}

Here we recall the minor summation formula \cite{IW1}.

\begin{lemma} \label{lem:minor-sum}
Let $X$ be a $2n \times N$ matrix and $A$ be an $N \times N$
 skew-symmetric matrix.
Then we have
$$
\sum_I \Pf \Delta^I_I(A) \det \Delta_I(X) = \Pf ( X A \trans X),
$$
where $I$ runs over all $2n$-element subsets of $[N]$.
\end{lemma}

By applying this minor-summation formula, we obtain

\begin{prop} \label{prop:schur}
Let $B = (b_{ij})_{i,j \ge 0}$ be the skew-symmetric matrix,
 whose entries $b_{ij}$ are given in Lemma~\ref{lem:coeff}.
Then, for a partition $\lambda$ of length $\le 2n$, we have
\begin{equation}
\sum_{\substack{\mu \subset \square(n,e) \\ \nu \subset \square(n,f)}}
 \LR^\lambda_{\mu,\nu}
  s_{\mu^\dagger(n,e)}(\vectz) s_{\nu^\dagger(n,f)}(\vectw)
=
\Pf \Delta^{I(\lambda)}_{I(\lambda)}(B).
\label{pf-schur3}
\end{equation}
\end{prop}

\begin{demo}{Proof}
Apply Lemma~\ref{lem:minor-sum} to the matrix
 $X = \left( x_i^k \right)_{1 \le i \le 2n, k \ge 0}$
 and the skew-symmetric matrix $B$.
Since $\det \Delta_{I(\lambda)}(X) / \Delta(\vectx) = s_\lambda(\vectx)$,
 the left hand side of (\ref{pf-schur}) becomes
$$
\frac{1}{\Delta(\vectx)}
\Pf \left(
 (x_j - x_i)
 h_{e+n-1}(x_i,x_j,\vectz)
 h_{f+n-1}(x_i,x_j,\vectw)
\right)_{1 \le i, j \le 2n}
 =
\sum_\lambda \Pf \Delta^{I(\lambda)}_{I(\lambda)} (B) s_\lambda(\vectx).
$$
where $\lambda$ runs over all partitions of length $\le 2n$.
Here we note that $\Pf \Delta^{I(\lambda)}_{I(\lambda)}(B) = 0$
 unless $\lambda \subset \square(2n,e+f)$.

On the other hand, the right hand side of (\ref{pf-schur})
 is expanded in the Schur function basis $\{ s_\lambda(\vectx) \}$
 as follows.
It follows from (\ref{LR1}) that
\begin{align*}
s_{\square(n,e)}(\vectx, \vectz)
 &=
 \sum_{\mu \subset \square(n,e)}
  s_\mu(\vectx) s_{\mu^\dagger(n,e)}(\vectz),
\\
s_{\square(n,f)}(\vectx, \vectw)
 &=
 \sum_{\nu \subset \square(n,f)}
  s_\nu(\vectx) s_{\nu^\dagger(n,f)}(\vectw).
\end{align*}
Hence we see that the right hand side of (\ref{pf-schur2}) becomes
\begin{align*}
s_{\square(n,e)}(\vectx, \vectz)
s_{\square(n,f)}(\vectx, \vectw)
 &=
\sum_{\substack{\mu \subset \square(n,e) \\ \nu \subset \square(n,f)}}
 s_\mu(\vectx) s_\nu(\vectx)
 s_{\mu^\dagger(n,e)}(\vectz) s_{\nu^\dagger(n,f)}(\vectw)
\\
&=
\sum_\lambda \left(
 \sum_{\substack{\mu \subset \square(n,e) \\ \nu \subset \square(n,f)}}
 \LR^\lambda_{\mu,\nu}
  s_{\mu^\dagger(n,e)}(\vectz) s_{\nu^\dagger(n,f)}(\vectw)
 \right)
 s_\lambda(\vectx).
\end{align*}

Comparing the coefficient of $s_\lambda(\vectx)$ on both sides
 of (\ref{pf-schur2}) completes the proof of (\ref{pf-schur3}).
\end{demo}

Now we can finish the proof of Theorem~\ref{thm:LR}.

\begin{demo}{Proof of Theorem~\ref{thm:LR}}
In the above argument, we take $p \ge n$ and $r = 0$.
In this case, the variables $\vectw$ disappear and we see that
$$
b_{kl} = \begin{cases}
 h_{(e+n-1)+(f+n-1)+1-k-l}(\vectz)
 &\text{if $0 \le k \le \min (e+n-1,f+n-1)$ and $l \ge f+n-1$}, \\
 0 &\text{otherwise}
\end{cases}
$$
and the equation (\ref{pf-schur3}) becomes
$$
\sum_{\mu \subset \square(n,e)}
 \LR^\lambda_{\mu,\square(n,f)} s_{\mu^\dagger(n,e)}(\vectz)
=
\Pf \Delta^{I(\lambda)}_{I(\lambda)}(B).
$$

The skew-symmetric matrix $B$ has the form
$$
B
 = \begin{pmatrix} O & C & O \\ -\trans C & O & O \\ O & O & O \end{pmatrix},
\quad
C
 = \left( h_{e+n-1-i-j}(\vectz) \right)_{0 \le i \le f+n-1, 0 \le j \le e+n-1}.
$$
>From the relation (\ref{pf-det}), we see that
 the subpfaffian $\Pf \Delta^{I(\lambda)}_{I(\lambda)}(B)$ vanishes unless
$$
\lambda_{n+1} + n-1 \le \min (e+n-1,f+n-1),
\quad
\lambda_n + n \ge f + n,
$$
i.e.,
$$
\lambda_{n+1} \le \min(e,f),
\quad
\lambda_n \ge f.
$$
If these conditions are satisfied, then we have
\begin{align*}
\Pf \Delta^{I(\lambda)}_{I(\lambda)}(B)
 &=
(-1)^{n(n-1)/2}
 \det \left(
  h_{\beta_i - \alpha_{n+1-j} - i + (n+1-j)}(\vectz)
 \right)_{1 \le i, j \le n}
\\
 &=
 (-1)^{n(n-1)/2} (-1)^{n(n-1)/2}
 \det \left(
   h_{\beta_i - \alpha_j - i + j}(\vectz)
 \right)_{1 \le i, j \le n}
\\
 &= s_{\beta/\alpha}(\vectz).
\end{align*}
Hence we have
$$
\sum_{\mu \subset \square(n,e)}
 \LR^\lambda_{\mu,\square(n,f)} s_{\mu^\dagger(n,e)}(\vectz)
 = s_{\beta/\alpha}(\vectz).
$$
Comparing the coefficients of $s_{\mu^\dagger(n,e)}(\vectz)$
 completes the proof.
\end{demo}


\end{document}